# NONSEMIMARTINGALES: STOCHASTIC DIFFERENTIAL EQUATIONS AND WEAK DIRICHLET PROCESSES

By Rosanna Coviello and Francesco Russo

*Scuola Normale Superiore di Pisa and Université Paris 13,
and Université Paris 13*

In this paper we discuss existence and uniqueness for a one-dimensional time inhomogeneous stochastic differential equation directed by an $\mathbb{F}$-semimartingale $M$ and a finite cubic variation process $\xi$ which has the structure $Q + R$, where $Q$ is a finite quadratic variation process and $R$ is *strongly predictable* in some technical sense: that condition implies, in particular, that $R$ is *weak Dirichlet*, and it is fulfilled, for instance, when $R$ is independent of $M$. The method is based on a transformation which reduces the *diffusion* coefficient multiplying $\xi$ to 1. We use generalized Itô and Itô–Wentzell type formulae. A similar method allows us to discuss existence and uniqueness theorem when $\xi$ is a Hölder continuous process and $\sigma$ is only Hölder in space. Using an Itô formula for *reversible* semimartingales, we also show existence of a solution when $\xi$ is a Brownian motion and $\sigma$ is only continuous.

## 1. Introduction.

This paper deals with the study of stochastic differential equations driven by a process which is not a semimartingale. We aim at illustrating how, using different types of Itô or Itô–Wentzell formulae, it is possible to establish existence and uniqueness results for a stochastic differential equation driven by a nonsemimartingale $\xi$ with a multiplication factor $\sigma$. When the paths of $\xi$ have very few regularity, more regularity on $\sigma$ is required. On the contrary, if the Hölder regularity of $\xi$ is $\gamma > \frac{1}{2}$, $\sigma$ only needs to fulfill a Hölder regularity.

As we said, one of the achievements of the paper is constituted by an Itô–Wentzell formula for processes having a finite cubic variation. There are today an incredible amount of generalized Itô formulae and it would be for us almost impossible to quote them all. The standard situation can be found in [11] and [27]; see also [28]. Given a finite quadratic variation process $\xi$, and $f \in C^{1,2}([0,1] \times \mathbb{R})$, one expands $f(t, \xi_t)$ as follows:

$$(1) \qquad f(t, \xi_t) = f(0, \xi_0) + \int_0^t \partial_s f(s, \xi_s) ds + \int_0^t \partial_x f(s, \xi_s) d^\circ \xi_s,$$

where the integral with respect to $\xi$ is a symmetric integral; see Definition 2.6. In the literature there are generalizations in several directions, among them the following:

1. The case that $\xi$ is not of finite quadratic variation, for instance, $\xi$ is a finite cubic variation and $f$ of class $C^{1,3}$ (see, e.g., [7]) or $\xi$ is a fractional Brownian motion with Hurst index $H > \frac{1}{6}$, and $f$ is of class $C^6$; see, for example, [1, 14];
2. The case when $\xi$ is a (*reversible*) semimartingale, so essentially a classical process but $f$ is of class $C^1$; see, in general, [12, 26].

Itô formula for finite quadratic variation processes admits extensions of Itô–Wentzell type, as in [10], where the dependence in time is of semimartingale type. More precisely, it is possible to expand the process $X_t(\xi_t)$, where $X_t(x)$ is a family of semimartingales depending on a parameter with respect to a given filtration $\mathbb{F} = (\mathcal{F}_t)$, if for every fixed parameter $x$, the semimartingale $X_t(x)$ admits a representation as a classical stochastic integral with respect to some vector of driving $\mathbb{F}$-semimartingales $(N^1, \ldots, N^n)$, $\xi$ is $\mathbb{F}$-adapted, and the vector $(\xi, N^1, \ldots, N^n)$ has all its mutual brackets; see Definition 2.3. We generalize this result, establishing an Itô–Wentzell formula for a finite cubic variation process $\xi$, provided that some technical assumption on $(\xi, N^1, \ldots, N^n)$ is fulfilled, see hypothesis ($\mathcal{D}$) in Definition 3.6: we assume the existence of a filtration $\mathbb{H} \supseteq \mathbb{F}$, with respect to which the vector $(N^1, \ldots, N^n)$ is still a vector of semimartingales, such that $\xi$ is decomposable into the sum of two $\mathbb{H}$-adapted processes $Q$ and $R$, where $(Q, N^1, \ldots, N^n)$ has all its mutual brackets, and $R$ is *strongly predictable* with respect to $\mathbb{H}$; see Definition 3.5. In particular, $R$ is an $\mathbb{H}$-*weak Dirichlet* process in the sense of [7]. We recall that an $\mathbb{H}$-*weak Dirichlet* process is the sum of a continuous $\mathbb{H}$-local martingale and of an $\mathbb{H}$-adapted process $Q$ such that $[Q, N] = 0$ for every $\mathbb{H}$-semimartingale $N$. Recent developments on that subject appeared in [13] and [2]. The mentioned hypothesis on $R$ is verified in the following cases:

- $R$ is $\mathcal{F}_0$ measurable;
- $R$ is independent from $(N^1, \ldots, N^n)$ and the filtration generated by $(N^1, \ldots, N^n)$ and the whole process $R$ contains $\mathbb{F}$.

Among others, the calculus developed to perform the Itô–Wentzell formula helps us to clarify the structure of $\mathbb{F}$-*weak Dirichlet* processes if $\mathbb{F}$ is the natural filtration associated with a Brownian motion $W$. If $Q$ is an $\mathbb{F}$-adapted process and $[Q, W]$ has all its mutual brackets, the covariation $[Q, L]$ can be computed explicitly for every continuous $\mathbb{F}$-semimartingale $L$; see Proposition 3.9. This allows us to prove that a process $A$ is $\mathbb{F}$-*weak Dirichlet* if and only if it is the sum of an $\mathbb{F}$-local martingale and of an $\mathbb{F}$-adapted process $Q$, with $[Q, W] = 0$.

On the other hand, a stochastic differential equation of the form

$$(2) \qquad d^\circ X_t = \sigma(t, X_t)[d^\circ \xi_t + \beta(t, X_t)\, d^\circ M_t + \alpha(t, \xi_t)\, dV_t]$$

is considered, where $M$ is a local martingale, $V$ a bounded variation process, and $\xi$ is a finite cubic variation process with $(\xi, M)$ verifying hypothesis $(\mathcal{D})$. We show, in different cases, how it is possible to apply the Itô formula to reduce the *diffusion* coefficient $\sigma$ to 1, and to formulate existence and uniqueness of equation (2) by studying equations where the process $\xi$ appears only as an additive term. The improper terminology of diffusion coefficient will be indeed used in the whole paper. A particular case of that equation was considered in [7] when $\beta = 0$. There $\sigma$ was of class $C^3$, and the notion of solution for a process $X$ was somehow unnatural since it required that the couple $(X, \xi)$ was a symmetric vector Itô process. In the case $\sigma$ is bounded from below by a positive constant, that equation can be investigated with our techniques, weakening the assumptions on the coefficients, enlarging the class of uniqueness and improving the sense of solution avoiding the notion of symmetric vector Itô process.

In the literature stochastic differential equations of forward type as

$$(3) \qquad d^- X_t = \sigma(X_t)\, d^- \xi_t + \beta(t, X_t)\, dL_t$$

were solved operating via classical transformations, in the case $\xi$ has finite quadratic variation; see [25] for definition of forward integral. In [27] a first attempt was done when $L$ has bounded variation. Similar independent results were established in [32]. In [10] existence and uniqueness were studied in a class of processes $(X(t, \xi_t))$, where $X(t, x)$ is a family of semimartingale depending on a parameter and $L$ is a semimartingale. There the regularity of $\sigma$ was of $C^4$ type with $\sigma', \sigma''$ being bounded. In that framework our result enlarges again the class of uniqueness, and we also require less regularity.

Equations looking similar to (2) were considered in the framework of T. Lyons and collaborators rough paths theory (see [21]) even in the multidimensional case when $\sigma$ is Lipschitz, $\alpha = 0$, for a process with deterministic $p$-variation strictly smaller than 3. Interesting reformulations of that integration theory and calculus with some applications to SDEs are given in [16] and [8]. Rough path analysis is purely deterministic in contrast with ours which couples the *pathwise* techniques of stochastic calculus via regularization and probabilistic concepts; see hypothesis $(\mathcal{D})$.

Another topic of interest is the study of equation

$$(4) \qquad d^\circ X_t = \sigma(t, X_t)[d^\circ \xi_t + \alpha(t, X_t)\, dt],$$

where $\sigma$ is only locally Hölder continuous, $\alpha$ is locally Lipschitz with linear growth, and $\xi$ is a $\gamma$-Hölder continuous process with $\gamma > \frac{1}{2}$. We apply the same method to this equation exploiting an Itô formula available for processes having Hölder continuous paths established in [31]. To this extent, we need to show that the symmetric integral of a process $f$ with respect to a process $g$ being Hölder continuous, respectively, of order $\gamma$ and $\delta$, with $\gamma + \delta > 1$, is the type of integral studied in [31]. Indeed, we prove that this integral is a particular case of the so-called *Young* integral introduced in a more general setting in [30]. Since the trajectories of the fractional Brownian motion are $\gamma$-Hölder continuous for every $\gamma$ strictly smaller than the Hurst parameter $H$, we are naturally induced to treat equations driven by the fractional Brownian motion with Hurst parameter $H > \frac{1}{2}$. Moreover, we combine our method with a recent result obtained in [22] with respect to an equation driven by a fractional Brownian motion with diffusion coefficient equal to 1. This permits us to improve our general result about existence and uniqueness of equation (2) when $\xi = B^H$, and $B^H$ is a fractional Brownian motion with Hurst index bigger than $\frac{1}{2}$, that is,

$$(5) \qquad d^\circ X_t = \sigma(t, X_t)[d^\circ B_t^H + \alpha(t, X_t)\, dt].$$

If the fractional Brownian motion reduces to a Brownian motion ($H = \frac{1}{2}$), an Itô formula for $C^1$ functions of *reversible* semimartingales is taken into consideration to formulate an existence theorem for equation (5), when $\sigma$ is only continuous and $\alpha$ is bounded measurable.

If $H$ is smaller than $\frac{1}{2}$, the Itô formula for Young type integral is no longer available. In spite of this, starting from our analysis, conditions to insure existence and uniqueness for equation (5) can still be deduced, treating the fractional Brownian motion with Hurst parameter $H \geq \frac{1}{3}$ as a strong finite cubic variation process. Essentially, in this case, the coefficient $\sigma$ is required to admit second continuous derivative with respect to the space variable.

On the other hand, remaining in the pure pathwise spirit, the Hölder nature of the fractional Brownian motion can be exploited to study equations of type (5), even when the Hurst parameter $H$ is smaller than $\frac{1}{2}$. The natural prolongation of Young integration and calculus is indeed rough path analysis.

Recently, several efforts were made in this direction (see [3, 4, 16]) to adapt results on rough paths theory to stochastic differential equations driven either by Hölder continuous processes with parameter $\gamma > \frac{1}{3}$ or by fractional Brownian motion with Hurst index $H > \frac{1}{4}$.

In [16] the author investigates existence and uniqueness of differential equations of type (4) with $\alpha = 0$, driven by irregular paths with Hölder exponent $\gamma$ greater than $\frac{1}{3}$. The multiplicative nonlinearity $\sigma$ was required to be differentiable till order two with second derivative $\delta$-Hölder continuous with $\delta > \frac{1}{\gamma} - 2$. At our knowledge, the first attempts to apply rough paths theory to the study of a stochastic differential equation driven by fractional Brownian motion of type (5) with $H$ strictly smaller than $\frac{1}{2}$ is constituted by [4]. There the authors considered the case $\frac{1}{4} < H < \frac{1}{2}$, and again $\alpha = 0$. They presented a pathwise approach to the solution of stochastic differential equations based on the so-called *universal limit theorem* established in [20]. To apply that result, the multiplicative coefficient $\sigma$ was assumed to be differentiable with bounded derivatives till order $[\frac{1}{H}] + 1$.

In both of the above-mentioned papers stochastic differential equations are solved in some specific setting and it is not obvious to see which kind of stochastic integral is involved.

A first result offering a link between the deterministic approach and the stochastic one can be found in [3]. There equation (5) is considered with $\alpha$ and $\sigma$ time independent vector fields. Assuming $\sigma$ differentiable and bounded till order $[\delta]$ with its $[\delta]$-derivative $(\delta - [\delta])$-Hölder, for some $\delta > \frac{1}{H}$, it is proved that the unique solution originated by the rough path method is actually a solution in some Stratonovich sense.

We come back to our paper. Our analysis of uniqueness, in the case of weak assumption on the *diffusion* coefficient, is inspired by classical ordinary differential equations of the type

$$(6) \qquad \frac{dX(t)}{dt} = \sigma(X(t)),$$

with $\sigma$ only continuous with linear growth. In the that case, Peano theorem insures existence but not uniqueness. Suppose that $\{x_0\} = \{x \in \mathbb{R}, \text{ s.t. } \sigma(x) = 0\}$. Then, if for some $\varepsilon > 0$,

$$(7) \qquad \int_{x_0}^{x_0+\varepsilon} \frac{1}{|\sigma|}(y)\, dy = \int_{x_0-\varepsilon}^{x_0} \frac{1}{|\sigma|}(y)\, dy = +\infty,$$

for every initial condition, this equation admits a unique solution. If the previous condition is not verified, then it is possible to show that at least two solutions for equation (6) exist, with initial condition $X_0 = x_0$. Suppose, for instance, that the second integral is finite. Setting $H(x) = \int_{x_0}^{x} \frac{1}{\sigma(y)}\, dy$, $x > x_0$, one can construct two solutions, that is, $X(t) \equiv x_0$ and $X(t) = H^{-1}(t)$. This phenomenon will be illustrated in the stochastic case, even with $\sigma$ inhomogeneous; see, for instance, Proposition 4.30 and Remark 4.31.

We observe that a similar condition as (7) appears in the study of one-dimensional stochastic differential equation of Itô type $dX(t) = \sigma(X(t))\, dW(t)$, where $W$ is a classical Brownian motion. Uniqueness for every initial condition holds if and only if

$$(8) \qquad \int_{x_0-\varepsilon}^{x_0+\varepsilon} \frac{1}{\sigma^2}(t)\, dt = +\infty,$$

for every $x_0 \in \mathbb{R}$; see [6].

To summarize, toward the study of equation (2), we innovate along the following axes with respect to the literature:

- We suppose that $\xi$ is a finite cubic variation process and $\sigma$ is time inhomogeneous.
- The notion of solution is clarified and we do not need to introduce the notion of symmetric vector Itô process.
- One new tool that we establish is an Itô–Wentzell type formula where finite cubic variation processes are involved.
- We continue the analysis related to the structure of weak Dirichlet processes.
- When the paths of $\xi$ are Hölder, with parameter greater than $\frac{1}{2}$, we require very weak regularity on the coefficients.
- In the case of classical Brownian motion, a new existence theorem is established for the Stratonovich equation.
- We drastically weaken the classical assumptions on the coefficients for existence and uniqueness. Our regularity assumptions are generally weaker than those intervening in rough path theory.

The paper is organized as follows. In Section 2 we recall some definitions and results about stochastic calculus with respect to finite cubic variation processes. We state the Itô formula and a result of stability of finite cubic variation through $C^1$ transformations. We also show some technical properties of the symmetric integral regarding its behavior when it is restricted to some subspace of the reference probability space, *stopped* or *shifted* with respect to some random time.

Section 3 deals with the class $\mathcal{C}_\xi^k(\mathbb{H})$ of the processes $Z$ so defined

$$Z_t = X_t(\xi_t),$$

being $X_t(x)$ an Itô field driven by a vector $(N^1, \ldots, N^n)$ of semimartingales such that hypothesis $(\mathcal{D})$ is verified for $(\xi, N^1, \ldots, N^n)$ (see Definition 3.1), with respect to some filtration $\mathbb{H}$, with regularity of order $k$ in the space variable. We prove that, if $\xi$ has a finite cubic variation, processes in $\mathcal{C}_\xi^1(\mathbb{H})$ still have finite cubic variation, and it is possible to establish an Itô–Wentzell formula to expand processes in $\mathcal{C}_\xi^3(\mathbb{H})$. In this section we also discuss connections with *weak Dirichlet* processes. We conclude this part proving the existence of the symmetric integral of a process in $\mathcal{C}_\xi^2(\mathbb{H})$ with respect to a process in $\mathcal{C}_\xi^2$, and using this result to formulate a chain-rule formula.

Section 4 discusses uniqueness and existence of equation (2). It is divided into nine subsections. The first and the second parts specify the notion of solution and describe the framework: we restrict ourselves to the case where the support $S$ of $\sigma$ is time-independent and a nonintegrability condition around its zeros of type (7) is fulfilled. The third part focuses on trajectories of solutions: if $X$ is a solution of equation (2), it can be expressed as a function of $\xi$ and a semimartingale. Moreover, its trajectories are forced to live in some connected component of $S$, as soon as the initial condition does. In the case the coefficients driving the equation are autonomous, a solution starting in $D = \mathbb{R}/S$ is identically equal to the initial condition. Putting things together, in the fourth part, we establish an equivalence between equation (2) and an equation of the same form but with diffusion coefficient equal to 1. We finally give some conditions for existence and uniqueness of this last equation. In the fifth subsection we use results of Section 3 to show that, under additional assumptions on the regularity of $\sigma$ and $\beta$, equation (2) admits a unique integral solution in the set $\mathcal{C}^2_\xi$. In the sixth one we revisit our results in the case $\xi$ has finite quadratic variation, and the symmetric integral is substituted by the forward integral. The seventh subsection is devoted to the application of the method when processes have Hölder trajectories. Subsection eight describes how it is possible to combine the result of [22] and ours to treat the specific case of an equation driven by fractional Brownian motion. Finally we discuss existence of solutions for a Stratonovich equation driven by a Brownian motion, with continuous diffusion coefficient and bounded measurable drift.

**2. Definitions, notation and basic calculus.** In this section we recall basic concepts and results about calculus with respect to finite cubic variation processes which will be useful later. For a more complete description of these arguments, the reader may refer to [7] or [14]. Throughout the paper $(\Omega, \mathcal{F}, P)$ will be a fixed probability space. All processes are supposed to be continuous and indexed by the time variable $t$ in $[0, 1]$. We adopt the notation $X_t = X_{(t \vee 0) \wedge 1}$, for every $t$ in $\mathbb{R}$. A sequence of continuous processes $(X^\varepsilon)_{\varepsilon > 0}$ will be said to converge *ucp* (*uniformly convergence in probability*) to a process $X$, if $\sup_{0 \le t \le 1} |X^\varepsilon_t - X_t|$ converges to zero in probability, when $\varepsilon$ goes to 0.

In the paper $C^{h,k}$ will be the space of all continuous functions $f : [0, 1] \times \mathbb{R} \to \mathbb{R}$, which are of class $C^h$ in $t$, with derivatives in $t$ up to order $h$ continuous in $(t, x)$, and of class $C^h$ in $x$, with derivatives in $x$ up to order $k$ continuous in $(t, x)$.

Let $n \ge 2$ and $(X^1, \ldots, X^n)$ be a vector of continuous processes. For any $\varepsilon > 0$ and $t$ in $[0, 1]$, set

$$[X^1, X^2, \ldots, X^n]_\varepsilon(t) = \frac{1}{\varepsilon} \int_0^t \prod_{k=1}^n (X^k_{s+\varepsilon} - X^k_s) \, ds$$

and

$$\|[X^1, X^2, \ldots, X^n]\|_\varepsilon = \frac{1}{\varepsilon} \int_0^1 \prod_{k=1}^n |X^k_{s+\varepsilon} - X^k_s| \, ds.$$

If $[X^1, X^2, \ldots, X^n]_\varepsilon(t)$ converges $ucp$, when $\varepsilon \to 0$, then the limiting process is called the *n-covariation* process of the vector $(X^1, \ldots, X^n)$ and denoted $[X^1, X^2, \ldots, X^n]$. If, furthermore, every subsequence $(\varepsilon_k)_{k \geq 0}$ admits a subsequence $(\bar\varepsilon_k)_{k \geq 0}$ such that

$$(9) \qquad \sup_{k \geq 0} \|[X^1, X^2, \ldots, X^n]\|_{\bar\varepsilon_k} < +\infty \qquad \text{a.s.},$$

then the *n*-covariation is said to exist in the *strong sense*. If the processes $(X^k)_{k=1}^n$ are all equal to a real valued process $X$, then the *n*-covariation of the considered vector will be denoted by $[X; n]$ and called the *n-variation* process. If $n = 2$, this process is the *quadratic variation* and it is denoted by $[X]$ or $[X, X]$. If $n = 3$, we will speak about *cubic variation*. If $X$ has a quadratic (resp., strong cubic) variation, $X$ will be called *finite quadratic variation* (resp., *strong cubic variation*) process.

REMARK 2.1. In [7] a different version of the definition of the strong *n*-variation is given. However, results contained there and recalled in the sequel can be proved to hold even under our weaker assumption.

EXAMPLE 2.2. We present several examples of strong finite cubic variation processes:

1. Let $(B_t^H, 0 \leq t \leq 1)$ be a fractional Brownian motion of Hurst index $H$, that is, a Gaussian process with zero mean and covariance

$$\text{Cov}(B_s^H, B_t^H) = \tfrac{1}{2}(s^{2H} + t^{2H} - |t - s|^{2H}).$$

   It follows from Remark 2.8 of [7] that the fractional Brownian motion with Hurst parameter $H = \tfrac{1}{3}$ is a strong cubic variation process.

2. Let $(B_t^{H,K}, 0 \leq t \leq 1)$ be a bifractional Brownian motion with parameters $H \in ]0, 1[, K \in ]0, 1]$. We recall (see [17]) that $B^{H,K}$ is a Gaussian process with zero mean and covariance

$$R(t, s) = \frac{1}{2^K}((t^{2H} + s^{2H})^K - |t - s|^{2HK}).$$

   In [29] is shown that $B^{H,K}$ is a strong finite cubic variation process if $HK \geq \tfrac{1}{3}$.

3. Let $(X_t, 0 \leq t \leq 1)$ be a Gaussian mean zero process starting at zero, with stationary increments. Set $V(t)^2 := \text{Var}(X_t)$, for every $t$ in $[0, 1]$. The Fubini theorem and the fact that the increments of $X$ are stationary permit us to perform the following evaluation:

$$E[\|X, X, X\|_\varepsilon] = \frac{c}{\varepsilon}(V(\varepsilon))^3,$$

   for some positive constant $c$. If, furthermore, $V(t) = O(t^{1/3})$, condition (9) holds. Moreover, using similar methods as in [15], it is possible to prove that the sequences of processes $[X, X, X]_\varepsilon$ converges $ucp$. In particular, $X$ is a strong cubic variation process.

4. Using [7], it is possible to exhibit examples of non-Gaussian strong finite cubic variation processes. One such process is of the type $X_t = \int_0^t G(t, s) \, dM_s$, where $M$ is a local martingale and $G$ is a continuous random field independent from $M$, essentially such that $[G(\cdot, s_1), G(\cdot, s_2), G(\cdot, s_3)]$ exist for any $s_1, s_2, s_3$. For example, one may choose $G(t, s) = B_{t-s}^H$, where $B^H$ is a fractional Brownian motion independent of $M$, with $H \geq \frac{1}{3}$.

DEFINITION 2.3. A vector $(X^1, X^2, \ldots, X^m)$ of continuous processes is said to have all its *mutual* (resp., *strong*) *n-covariations* if $[X^{i_1}, X^{i_2}, \ldots, X^{i_n}]$ exists (resp., exists in the strong sense) for any choice (even with repetition) of indices $i_1, i_2, \ldots, i_n$ in $\{1, 2, \ldots, m\}$. If $n = 2$, we will also say that the vector $(X^1, X^2, \ldots, X^m)$ has all its *mutual brackets*. In that case $[X^1, \ldots, X^m]$ has bounded variation.

PROPOSITION 2.4. *If condition* (9) *holds, then* $[X^1, X^2, \ldots, X^n]$ *has bounded variation whenever it exists.*

REMARK 2.5. 1. If the $n$-variation $[X; n]$ exists in the strong sense for some $n$, then $[X; m] = 0$ for all $m > n$. In particular, since the 2-covariation of two semimartingales exists strongly and agrees with their usual covariation (see [25]), for any semimartingale $S$, $[S; n] = 0$ for all $n \geq 3$.

2. Let $(X^1, \ldots, X^n)$ be a vector having a strong $n$-covariation, and $Y$ a continuous process. Then

$$\frac{1}{\varepsilon} \int_0^{\cdot} Y_s \prod_{k=1}^n (X_{s+\varepsilon}^k - X_s^k) \, ds$$

converges $ucp$ to

$$\int_0^{\cdot} Y \, d[X^1, X^2, \ldots, X^n].$$

3. If $(X^1, \ldots, X^n)$ has its strong $n$-covariation then for every vector of continuous processes $(Y^1, Y^2, \ldots, Y^m)$, the vector

$$(X^1, \ldots, X^n, Y^1, \ldots, Y^m)$$

has its strong $(n + m)$-covariation equal to zero.

4. If the $n$-variation $[X; n]$ exists in the strong sense, then for every continuous process $Y$ and every $m > n$ such that $[Y; m]$ exists in the strong sense, we have

$$[X, \overbrace{Y, Y, \ldots, Y}^{(m-1)\text{times}}] = 0.$$

DEFINITION 2.6. Let $X$ and $Y$ be two continuous processes. For any $\varepsilon > 0$ and $t$ in $[0, 1]$, set

$$I_\varepsilon^\circ(t, X, Y) = \frac{1}{2\varepsilon} \int_0^t Y_s(X_{s+\varepsilon} - X_{s-\varepsilon})\, ds.$$

If the process $I_\varepsilon^\circ(\cdot, X, Y)$ converges $ucp$ when $\varepsilon$ goes to zero, then the limiting process will be denoted by $\int_0^t Y\, d^\circ X$ and called the *symmetric integral*.

REMARK 2.7. 1. It is easy to show that the symmetric integral, if it exists, is the limit $ucp$ of

$$J_\varepsilon^\circ(t) = \frac{1}{2\varepsilon} \int_0^t (Y_{s+\varepsilon} + Y_s)(X_{s+\varepsilon} - X_s)\, ds.$$

2. Let $X$ be a continuous semimartingale with respect to some filtration $\mathbb{F}$ and $Y$ an $\mathbb{F}$-adapted continuous process such that $[X, Y]$ exists. Then the symmetric integral $\int_0^\cdot Y_s\, d^\circ X_s$ exists,

$$\int_0^\cdot Y_s\, d^\circ X_s = \int_0^\cdot Y_s\, dX_s + \frac{1}{2}[X, Y],$$

and it coincides with classical Stratonovich integral if $Y$ is an $\mathbb{F}$-semimartingale.

We conclude this section by recalling a result about stability of the strong $n$-covariation through $C^1$ transformations, the Itô formula for *strong cubic variation* processes and a *chain-rule* formula, all of them established in [7], Propositions 2.7, 3.7 and Lemma 3.18.

PROPOSITION 2.8. *Let* $F^1, \ldots, F^n$ *be* $n$ *functions in* $C^1(\mathbb{R}^n)$. *Let* $X = (X^1, \ldots, X^n)$ *be a vector of continuous processes having all its mutual strong $n$-covariations. Then the vector*

$$(F^1(X), \ldots, F^n(X))$$

*has the same property and*

$$[F^1(X), \ldots, F^n(X)] = \sum_{1 \le i_1, \ldots, i_n \le n} \int_0^t \partial_{i_1} F^1(X) \cdots \partial_{i_n} F^n(X)\, d[X^{i_1}, \ldots, X^{i_n}].$$

PROPOSITION 2.9. *Let* $V = (V^1, \ldots, V^m)$ *be a vector of bounded variation processes and* $\xi$ *be a strong cubic variation process. Then for every* $F$ *belonging to the class* $C^{1,3}(\mathbb{R}^m \times \mathbb{R})$, *it holds*

$$F(V_t, \xi_t) = F(V_0, \xi_0) + \sum_{i=1}^m \int_0^t \partial_{V^i} F(V_s, \xi_s)\, dV_s^i + \int_0^t \partial_\xi F(V_s, \xi_s)\, d^\circ \xi_s$$

$$- \frac{1}{12} \int_0^t \partial_\xi^{(3)} F(V_s, \xi_s)\, d[\xi, \xi, \xi]_s.$$

LEMMA 2.10. *Let $\xi$ be a strong cubic variation process. Suppose that $\psi$ and $\phi$ are, respectively, in $C^{1,3}([0,1] \times \mathbb{R})$ and $C^{1,2}([0,1] \times \mathbb{R})$. Then $X = \int_0^{\cdot} \phi(s, \xi_s) \, d^{\circ}\xi_s$, and $\int_0^{\cdot} \psi(s, \xi_s) \, d^{\circ}X_s$ exist and*

$$\int_0^{\cdot} \psi(s, \xi_s) \, d^{\circ} X_s = \int_0^{\cdot} \phi\psi(s, \xi_s) \, d^{\circ}\xi_s - \frac{1}{4} \int_0^{\cdot} \partial_{\xi}\psi \, \partial_{\xi}\phi(s, \xi_s) \, d[\xi, \xi, \xi]_s.$$

In the sequel of the paper we will need to deal with the restriction of symmetric integrals to subspaces of $\Omega$, as well as with symmetric integrals *stopped* or *shifted* with respect to random times. We list some simple technical properties about these operations.

If $B$ is an element of $\mathcal{F}$, with $P(B) > 0$, $\mathcal{F}^B$ will denote the restriction of $\mathcal{F}$ on $B : \mathcal{F}^B = \{F \cap B, F \in \mathcal{F}\}$, $P^B$ the probability measure conditioned on $B$, and if $f$ is a random variable on $(\Omega, \mathcal{F}, P)$, $f^B$ will denote the restriction of $f$ to $B$.

LEMMA 2.11. *Let $B$ in $\mathcal{F}$ with $P(B) > 0$. Let $X$ and $Y$ be two continuous processes such that $\int_0^{\cdot} X \, d^{\circ}Y$ exists. Then $\int_0^{\cdot} X^B \, d^{\circ}Y^B$ exists and*

$$\int_0^{\cdot} X_t^B \, d^{\circ}Y_t^B = \left( \int_0^{\cdot} X_t \, d^{\circ}Y_t \right)^B, \qquad P^B \ a.s.$$

PROOF. The result follows immediately after having observed that for every $\delta > 0$,

$$P^B\left( \left\{ \sup_{t \in [0,1]} \left| I_{\varepsilon}^{\circ}(t, X^B, Y^B) - \left( \int_0^t X_s \, d^{\circ}Y_s \right)^B \right| > \delta \right\} \right)$$

$$\leq \frac{1}{P(B)} P\left( \left\{ \sup_{t \in [0,1]} \left| I_{\varepsilon}^{\circ}(t, X, Y) - \left( \int_0^t X_s \, d^{\circ}Y_s \right) \right| > \delta \right\} \right). \qquad \square$$

If $\tau$ and $X$ are, respectively, a random time and a stochastic process on $(\Omega, \mathcal{F}, P)$, $X^{\tau}$ will denote the stochastic process $X$ stopped to time $\tau : X_t^{\tau} = X_{t \wedge \tau}, 0 \leq t \leq 1$.

LEMMA 2.12. *Let $\tau$ be a random time on $(\Omega, \mathcal{F}, P)$, with $P(\tau \leq 1) = 1$, $X$ and $Y$ two continuous stochastic processes such that $\int_0^{\cdot} X \, d^{\circ}Y$ exists. Then it holds:*

$$\int_0^{\cdot} X_s^{\tau} \, d^{\circ}Y_s^{\tau} = \left( \int_0^{\cdot} X_s \, d^{\circ}Y_s \right)^{\tau};$$

$$\int_0^{\cdot} X_{\tau+s} d^{\circ}(Y_{\tau+s}) = \int_{\tau}^{\tau+\cdot} X_s \, d^{\circ}Y_s.$$

PROOF. We clearly have

$$\sup_{t \in [0,1]} \left| I_\varepsilon^\circ(t \wedge \tau, X, Y) - \left( \int_0^\cdot X_s \, d^\circ Y_s \right)_{t \wedge \tau} \right| \leq \sup_{t \in [0,1]} \left| I_\varepsilon^\circ(t, X, Y) - \int_0^t X_s \, d^\circ Y_s \right|.$$

Therefore, for the first part of the statement, we have to show that $\lim_{\varepsilon \to 0} a_\varepsilon = 0$, in probability, with

$$a_\varepsilon = \sup_{t \in [0,1]} |I_\varepsilon^\circ(t \wedge \tau, X, Y) - I_\varepsilon^\circ(t, X^\tau, Y^\tau)|.$$

We can write

$$a_\varepsilon \leq \sup_{t \in [0,1]} \left| \frac{1}{\varepsilon} \int_{(\tau - \varepsilon) \wedge t}^{\tau \wedge t} X_s (Y_\tau - Y_{s+\varepsilon}) \, ds \right|$$

$$+ \sup_{t \in [0,1]} \left| \frac{1}{\varepsilon} \int_{\tau \wedge t}^{(\tau + \varepsilon) \wedge t} X_\tau (Y_\tau - Y_{s-\varepsilon}) \, ds \right|.$$

The convergence to zero almost surely, and so in probability, of the sequence of processes $(a_\varepsilon)$ is due to the continuity of the processes $X$ and $Y$.

The second statement is a straightforward consequence of a simple change of variables which let us obtain $I^\circ(t, X_{\tau+\cdot}, Y_{\tau+\cdot}) = I^\circ(\tau + \cdot, X, Y) - I^\circ(\tau, X, Y)$. □

By similar arguments, it is also possible to show the following lemma.

LEMMA 2.13. *Let* $(X^1, \ldots, X^n)$ *be a vector of continuous processes having its $n$-covariation, $\tau$ a random time with $P(\tau \leq 1) = 1$, and $B$ an element of $\mathcal{F}$. Then the vectors* $((X^1)^B, \ldots, (X^n)^B)$, $((X^1)^\tau, \ldots, (X^n)^\tau)$ *and* $(X^1_{\tau+\cdot}, \ldots, X^n_{\tau+\cdot})$ *have their $n$-covariation and*

$$[X^1, \ldots, X^n]^B = [(X^1)^B, \ldots, (X^n)^B], \qquad P^B \, a.s.;$$

$$[X^1, \ldots, X^n]^\tau = [(X^1)^\tau, \ldots, (X^n)^\tau];$$

$$[X^1_{\tau+\cdot}, \ldots, X^n_{\tau+\cdot}] = [X^1, \ldots, X^n]_{\tau+\cdot} - [X^1, \ldots, X^n]_\tau.$$

**3. Itô-fields evaluated at strong cubic variation processes.**

3.1. *Stability of strong cubic variation.* At this stage we introduce some definitions adapted from [10], which treated the finite quadratic variation case. From now on $\mathbb{H} = (\mathcal{H}_t)_{t \in [0,1]}$ will denote a filtration on $(\Omega, \mathcal{F})$, satisfying the *usual assumptions*.

DEFINITION 3.1. A random field $(X(t,x), 0 \leq t \leq 1, x \in \mathbb{R})$ is called a $C^k$ $\mathbb{H}$-*Itô-martingale field* driven by the vector $N = (N^1, \ldots, N^n)$, if $N$ is a vector of local martingales with respect to $\mathbb{H}$, and

$$(10) \qquad X(t,x) = f(x) + \sum_{i=1}^{n} \int_0^t a^i(s,x)\, dN_s^i,$$

where

$f : \Omega \times \mathbb{R} \to \mathbb{R}$ is, for every $x$, $\mathcal{H}_0$-measurable and belonging to $C^k(\mathbb{R})$ $a.s$;
$X$ and $a^i : [0,1] \times \mathbb{R} \times \Omega \to \mathbb{R}$, $i = 1, \ldots, n$ are $\mathbb{H}$-adapted for every $x$, almost surely continuous with their partial derivatives with respect to $x$ in $(t,x)$ up to order $k$;
for every index $h \leq k$, it holds

$$\partial_x^{(h)} X(t,x) = \partial_x^{(h)} f(x) + \sum_{i=1}^{n} \int_0^t \partial_x^{(h)} a^i(s,x)\, dN_s^i.$$

DEFINITION 3.2. Let $p \geq 1$. A continuous random field $(Z(t,x), 0 \leq t \leq 1, x \in \mathbb{R})$, is called an $\mathbb{H}$-*strict zero $p$-variation process* if it is $\mathbb{H}$-adapted for every $x$, and

$$(11) \qquad \sup_{|x| \leq R} \frac{1}{\varepsilon} \int_0^1 |Z(t+\varepsilon, x) - Z(t,x)|^p \, dt \to 0 \qquad \text{in probability,}$$

for all $R > 0$.

If $p = 2$ (resp., $p = 3$), $Z$ will be called an $\mathbb{H}$-*strict zero quadratic* (resp., *cubic*) process.

Note that if

$$(12) \qquad Z(t,x) = \sum_{j=1}^{m} \int_0^t b^j(s,x)\, dV_s^j,$$

where $b^j$ are continuous fields, and $(V_t^j)_{0 \leq t \leq 1}$, $j = 1, \ldots, m$, are bounded variation processes, then (11) is verified for every $p > 1$.

DEFINITION 3.3. A random field $X$ will be called a $C^k$ $\mathbb{H}$-*Itô-semimartingale field* if it is the sum of a $C^k$ $\mathbb{H}$-Itô-martingale field and an $\mathbb{H}$-*strict zero quadratic variation* process $Z$ having the form (12):

$$(13) \qquad X(t,x) = f(x) + \sum_{i=1}^{n} \int_0^t a^i(s,x)\, dN_s^i + \sum_{j=1}^{m} \int_0^t b^j(s,x)\, dV_s^j,$$

with coefficients $(b^j)_{j=1}^{m}$ continuous with their partial derivatives with respect to $x$ in $(t,x)$ up to order $k$.

PROPOSITION 3.4. *Let $X = (X^i(t,x), 0 \le t \le 1, x \in \mathbb{R}, i = 1, 2, 3)$ be a vector of random fields being the sum of a vector of $C^1$ $\mathbb{H}$-Itô-martingale fields*

$$(Y^i(t,x), 0 \le t \le 1, x \in \mathbb{R}, i = 1, 2, 3),$$

*driven by the vector of local martingales $(N^1, \ldots, N^n)$, and of a vector of $\mathbb{H}$-strict zero cubic variation processes $(Z^i(t,x), 0 \le t \le 1, x \in \mathbb{R}, i = 1, 2, 3)$ which are a.s. in $C^{0,1}([0,1] \times \mathbb{R})$:*

$$X^i = Y^i + Z^i, \qquad i = 1, 2, 3.$$

*Let $\xi$ be a strong cubic variation and $\mathbb{H}$-adapted process. Then the vector $X(\cdot, \xi)$ has its strong mutual 3-covariations and*

$$[X^{i_1}(\cdot, \xi), X^{i_2}(\cdot, \xi), X^{i_3}(\cdot, \xi)] = \int_0^{\cdot} (\partial_x X^{i_1})(\partial_x X^{i_2})(\partial_x X^{i_3})(s, \xi_s) \, d[\xi, \xi, \xi]_s,$$

*for every choice of indices $(i_1, i_2, i_3)$ in $\{1, 2, 3\}$.*

PROOF. We first remark that it is not reductive to suppose that the vector of the driving local martingales is the same for all the Itô fields taken into consideration. We consider the case $X = X^1 = X^2 = X^3 = Y + Z$. The proof in the general case requires the same essential concepts. We suppose also, for simplicity of notation, that the $C^1$ $\mathbb{H}$-Itô-martingale field has the form (10) with $n = 1$, $N^1 = N$, $a^1 = a$. We have to prove that

$$C_\varepsilon = \frac{1}{\varepsilon} \int_0^{\cdot} (X(s + \varepsilon, \xi_{s+\varepsilon}) - X(s, \xi_s))^3 \, ds$$

converges *ucp* to $\int_0^{\cdot} (\partial_x X(s, \xi_s))^3 \, d[\xi, \xi, \xi]_s$, and that $X(\cdot, \xi)$ verifies condition (9). We can write

$$
\begin{aligned}
X(s + \varepsilon, \xi_{s+\varepsilon}) - X(s, \xi_s) &= (X(s + \varepsilon, \xi_{s+\varepsilon}) - X(s + \varepsilon, \xi_s)) \\
&\quad + (X(s + \varepsilon, \xi_s) - X(s, \xi_s)) \\
&= A(s, \varepsilon) + B(s, \varepsilon),
\end{aligned}
$$

so as to decompose $C_\varepsilon$ as follows:

$$C_\varepsilon(t) = I_\varepsilon^1(t) + I_\varepsilon^2(t) + 3I_\varepsilon^3(t) + 3I_\varepsilon^4(t),$$

with

$$I_\varepsilon^1(t) = \frac{1}{\varepsilon} \int_0^t (A(s, \varepsilon))^3 \, ds, \qquad\qquad I_\varepsilon^2(t) = \frac{1}{\varepsilon} \int_0^t (B(s, \varepsilon))^3 \, ds,$$

$$I_\varepsilon^3(t) = \frac{1}{\varepsilon} \int_0^t (A(s, \varepsilon))^2 (B(s, \varepsilon)) \, ds, \qquad I_\varepsilon^4(t) = \frac{1}{\varepsilon} \int_0^t (A(s, \varepsilon))(B(s, \varepsilon))^2 \, ds.$$

Since $X$ is differentiable in $\xi$, $A(s, \varepsilon)$ may be rewritten as

$$A(s, \varepsilon) = \rho(s, \varepsilon)(\xi_{s+\varepsilon} - \xi_s),$$

with

$$\rho(s, \varepsilon) = \int_0^1 \partial_x X(s + \varepsilon, \xi_s + \lambda(\xi_{s+\varepsilon} - \xi_s)) \, d\lambda.$$

Then

$$I_\varepsilon^1(t) = \frac{1}{\varepsilon} \int_0^t (\partial_x X(s, \xi_s))^3 (\xi_{s+\varepsilon} - \xi_s)^3 \, ds$$

$$+ \frac{1}{\varepsilon} \int_0^t ((\rho(s, \varepsilon))^3 - (\partial_x X(s, \xi_s))^3)(\xi_{s+\varepsilon} - \xi_s)^3 \, ds.$$

By Remark 2.5.2, the first term of this sum converges $ucp$ to

$$\int_0^\cdot (\partial_x X(s, \xi_s))^3 \, d[\xi, \xi, \xi]_s,$$

while the absolute value of the second term is bounded by

$$\sup_{s \in [0,1]} |(\rho(s, \varepsilon))^3 - (\partial_x X(s, \xi_s))^3| \left( \frac{1}{\varepsilon} \int_0^1 |\xi_{s+\varepsilon} - \xi_s|^3 \, ds \right),$$

which converges to zero in probability since $\partial_x X$ is continuous, and $\xi$ is a *strong cubic variation* process.

We show that $I_\varepsilon^2(t)$ converges to zero $ucp$. We observe that we can apply a substitution argument thanks to the Hölder continuity of $a$ (see [27], Proposition 2.1) and the adaptedness of the process $\xi$, and get

$$B(s, \varepsilon) = \left( \int_s^{s+\varepsilon} a(r, x) \, dN_r \right)_{x=\xi_s} + (Z(s + \varepsilon, \xi_s) - Z(s, \xi_s))$$

$$= \int_s^{s+\varepsilon} a(r, \xi_s) \, dN_r + (Z(s + \varepsilon, \xi_s) - Z(s, \xi_s)).$$

Then

$$|I_\varepsilon^2(t)| \leq \frac{1}{\varepsilon} \int_0^1 |B(s, \varepsilon)|^3 \, ds$$

$$\leq \frac{4}{\varepsilon} \int_0^1 \left| \int_s^{s+\varepsilon} a(r, \xi_s) \, dN_r \right|^3 \, ds$$

$$+ \frac{4}{\varepsilon} \int_0^1 |Z(s + \varepsilon, \xi_s) - Z(s, \xi_s)|^3 \, ds.$$

For every $k$ in $N^*$, we set

$$\Omega^k = \{[N]_1 \leq k\} \cap \left\{ \sup_{t \in [0,1]} |\xi_t| \leq k \right\},$$

$$\tau^k = \inf\{t | [N]_t \geq k\}, \qquad N^k = N^{\tau^k}.$$

Then $\tau^k$ is a stopping time and by optional sampling theorem, $N^k$ is a local square integrable martingale. Since $\bigcup_{k=0}^{\infty} \Omega^k = \Omega$, almost surely, it is sufficient to verify that, for every $k$ in $\mathbb{N}^*$, the sequence of processes $(I_{\Omega_k} I_\varepsilon^2(t))$ converges to zero $ucp$. Since $Z$ is an $\mathbb{H}$-strict zero cubic variation process and on $\Omega_k$ the process $\xi$ is bounded by a constant,

$$\lim_{\varepsilon \to 0} I_{\Omega^k} \left( \frac{1}{\varepsilon} \int_0^\cdot (Z(s+\varepsilon, \xi_s) - Z(s, \xi))^3 \, ds \right) = 0 \qquad ucp,$$

and so we get the desired convergence if

$$\lim_{\varepsilon \to 0} \int_0^1 \frac{1}{\varepsilon} \left| \int_s^{s+\varepsilon} a^k(r, \xi_s) \, dN_r^k \right|^3 ds = 0 \qquad \text{in probability},$$

where $a^k : [0,1] \times \mathbb{R} \to \mathbb{R}$ has the same regularity of $a$, it is bounded and it agrees with $a$ on $[0,1] \times \{x \in \mathbb{R} \mid |x| \le k\}$. We can write

$$\int_0^1 \frac{1}{\varepsilon} \left| \int_s^{s+\varepsilon} a^k(r, \xi_s) \, dN_r^k \right|^3 ds$$
$$\le \frac{4}{\varepsilon} \int_0^1 \left| \int_s^{s+\varepsilon} a^k(r, \xi_r) \, dN_r^k \right|^3 ds$$
$$+ \frac{4}{\varepsilon} \int_0^1 \left| \int_s^{s+\varepsilon} (a^k(r, \xi_s) - a^k(r, \xi_r)) \, dN_r^k \right|^3 ds.$$

The process $\int_0^\cdot a^k(r, \xi_r) \, dN_r^k$ is a continuous semimartingale, then it has a finite quadratic variation by Remark 2.5.1 and so the first term of the sum converges to zero in probability being bounded by

$$\left( \sup_{t \in [0,1]} \left| \int_s^{s+\varepsilon} a^k(r, \xi_r) \, dN_r^k \right| \right) \left( \int_0^1 \frac{1}{\varepsilon} \left| \int_s^{s+\varepsilon} a^k(r, \xi_r) \, dN_r^k \right|^2 ds \right).$$

Therefore, to conclude, we only need to apply the Burkholder inequality and the Lebesgue dominated convergence theorem to see that

$$\lim_{\varepsilon \to 0} E \left[ \int_0^1 \frac{1}{\varepsilon} \left| \int_s^{s+\varepsilon} (a^k(r, \xi_s) - a^k(r, \xi_r)) \, dN_r^k \right|^3 ds \right] = 0.$$

Finally, by the Hölder inequality,

$$|I_\varepsilon^3(t)| \le \left( \frac{1}{\varepsilon} \int_0^1 |A(s, \varepsilon)|^3 \, ds \right)^{2/3} \left( \frac{1}{\varepsilon} \int_0^1 |B(s, \varepsilon)|^3 \, ds \right)^{1/3}$$

and

$$|I_\varepsilon^4(t)| \le \left( \frac{1}{\varepsilon} \int_0^1 |A(s, \varepsilon)|^3 \, ds \right)^{1/3} \left( \frac{1}{\varepsilon} \int_0^1 |B(s, \varepsilon)|^3 \, ds \right)^{2/3},$$

then $I_\varepsilon^3(t)$ and $I_\varepsilon^4(t)$ converge to zero *ucp*, since, as already proved before, $\frac{1}{\varepsilon}\int_0^1 |B(s,\varepsilon)|^3\,ds$ converges to zero in probability and

$$(14) \qquad \frac{1}{\varepsilon}\int_0^1 |A(s,\varepsilon)|^3\,ds \leq \|\xi,\xi,\xi\|_\varepsilon \sup_{s\in[0,1]} |\rho(s,\varepsilon)|^3.$$

We conclude observing that the cubic variation of $X$ exists strongly thanks to inequality (14), the strong finite cubic variation of $\xi$ and the convergence to zero in probability of $\frac{1}{\varepsilon}\int_0^1 |B(s,\varepsilon)|^3\,ds$.  $\square$

3.2. *Strong predictability, covariations and weak Dirichlet processes.*  Given a vector of processes $(N^1,\dots,N^n)$, $\mathcal{S}(N^1,\dots,N^n)$ will denote the set of all filtrations on $(\Omega,\mathcal{F})$ with respect to which $(N^1,\dots,N^n)$ is a vector of semi-martingales.

DEFINITION 3.5. A process $R$ is *strongly predictable* with respect to $\mathbb{H}$ if

$$\exists\ \delta > 0, \text{ such that } R_{\varepsilon+}. \text{ is } \mathbb{H}\text{-adapted, for every } \varepsilon \leq \delta.$$

This notion constitutes in fact the direct generalization of the notion of predictability intervening in the discrete time case.

DEFINITION 3.6. We will say that the vector $(\xi, N^1,\dots,N^n)$ satisfies *hypothesis* ($\mathcal{D}$) with respect to $\mathbb{H}$, if $\mathbb{H}$ belongs to $\mathcal{S}(N^1,\dots,N^n)$, and there exist two continuous processes, adapted to $\mathbb{H}$, such that

$$\xi = R + Q;$$

($\mathcal{D}$)     $R$ is strongly predictable with respect to $\mathbb{H}$;

the vector $(Q, N^1,\dots,N^n)$ has all its mutual brackets.

We give two examples where there exists a filtration $\mathbb{H}$ with respect to which the decomposition ($\mathcal{D}$) occurs.

EXAMPLE 3.7. Let $(N^1,\dots,N^n)$ be a vector of local martingales with respect to a filtration $\mathbb{F} = (\mathcal{F}_t)_{t\in[0,1]}$. Suppose that $\xi = R + Q$, where

$R$ is $\mathbb{F}_0$-measurable;

$(Q, N^1,\dots,N^n)$ has all its mutual brackets and $Q$ is $\mathbb{F}$-adapted.

Then the hypothesis ($\mathcal{D}$) is satisfied with respect to the filtration $\mathbb{F}$.

EXAMPLE 3.8. Let $(N^1, \ldots, N^n)$ be a vector of semimartingales with respect to its natural filtration $\mathbb{G} = (\mathcal{G}_t)_{t \in [0,1]}$. Suppose that $\xi = R + Q$, where

> $R$ is independent from $(N^1, \ldots, N^n)$;
>
> $(Q, N^1, \ldots, N^n)$ has all its mutual brackets.

Then, if $Q$ is adapted to the filtration

$$\mathbb{H} = (\mathcal{G}_t \vee \sigma(R))_{t \in [0,1]},$$

the vector $(\xi, N^1, \ldots, N^n)$ satisfies the hypothesis ($\mathcal{D}$) with respect to $\mathbb{H}$.

For every $\mathbb{H}$-local martingale $N$, we denote with $\mathcal{L}_N^2(\mathbb{H})$ the set of all progressively measurable processes $h$ such that

$$\|h\|_{L^2(d[N])} = \int_0^1 h_s^2 \, d[N]_s < +\infty \qquad \text{a.s.}$$

$\mathcal{L}_N^2(\mathbb{H})$ endowed with the topology of the convergence in probability with respect to the norm $\|\cdot\|_{L^2(d[N])}$ is an $F$-space in the sense of [5]. The $F$-space of all continuous $\mathbb{H}$-adapted processes equipped with the uniform convergence in probability will be denoted by $\mathcal{A}(\mathbb{H})$.

PROPOSITION 3.9. *Let $Q$ be a continuous and $\mathbb{H}$-adapted process and $N$ a continuous $\mathbb{H}$-local martingale such that $(Q, N)$ has all its mutual brackets. Then for every $h$ in $\mathcal{L}_N^2(\mathbb{H})$, and $Y = \int_0^\cdot h_s \, dN_s$, the bracket $[Q, Y]$ exists and*

$$[Q, Y] = \int_0^\cdot h_s \, d[Q, N]_s.$$

*In particular, $(Q, Y)$ has all its mutual brackets and $[Q, Y]$ has bounded variation.*

PROOF OF PROPOSITION 3.9. By localization arguments, we do not lose generality if we suppose that $Q$ is uniformly bounded and $N$ is square integrable. We set $\Gamma(h) := \int_0^\cdot h_s \, d[Q, N]_s$, for every $h$ in $\mathcal{L}_N^2(\mathbb{H})$, and for every $\varepsilon > 0$, we consider the map $\Gamma_\varepsilon : \mathcal{L}_N^2(\mathbb{H}) \to \mathcal{A}(\mathbb{H})$ so-defined:

$$\Gamma_\varepsilon(h) = \frac{1}{\varepsilon} \int_0^\cdot (Q_{\varepsilon+s} - Q_s) \left( \int_s^{s+\varepsilon} h_r \, dN_r \right) ds.$$

$\Gamma_\varepsilon$ is a linear and continuous operator from $\mathcal{L}_N^2(\mathbb{H})$ to $\mathcal{A}(\mathbb{H})$. Let $h$ be continuous. We claim that $(\Gamma_\varepsilon(h))$ converges $ucp$ to $\Gamma(h)$. Remark 2.5.2 implies

$$\lim_{\varepsilon \to 0} \int_0^\cdot h_s (Q_{s+\varepsilon} - Q_s)(N_{s+\varepsilon} - N_s) = \Gamma(h) \qquad ucp.$$

We hence achieve the claim if

$$\lim_{\varepsilon \to 0} I^\varepsilon(t) = \lim_{\varepsilon \to 0} \left| \frac{1}{\varepsilon} \int_0^t (Q_{s+\varepsilon} - Q_s) \left( \int_s^{s+\varepsilon} (h_s - h_r) \, dN_r \right) ds \right| = 0 \qquad ucp.$$

Again by standard localization techniques, we can suppose $h$ uniformly bounded. We use the Cauchy–Schwarz inequality to write

$$I^\varepsilon(t) \leq \left( \frac{1}{\varepsilon} \int_0^\cdot (Q_{s+\varepsilon} - Q_s)^2 \, ds \right)^{1/2} \left( \int_0^\cdot \frac{1}{\varepsilon} \left| \int_s^{s+\varepsilon} (h_s - h_r) \, dN_r \right|^2 ds \right)^{1/2}.$$

The expectation of the second factor of the product is convergent to zero by the Burkholder inequality, the continuity and the boundness of $h$.

Moreover, it is possible to show that, for every $h$ in $\mathcal{L}_N^2(\mathbb{H})$,

$$\sup_{\varepsilon > 0} d_1(\Gamma_\varepsilon(h), 0) \leq d_2(h, 0),$$

being $d_1$ and $d_2$ two metrics inducing the given topologies of $\mathcal{A}(\mathbb{H})$ and $\mathcal{L}_N^2(\mathbb{H})$, respectively. We recall that $\mathbb{H}$-adapted continuous processes are dense in $\mathcal{L}_N^2(\mathbb{H})$, so that the Banach–Steinhaus theorem for Fréchet spaces ([5], Chapter 2.1) and the density of continuous processes permit us to conclude. $\square$

PROPOSITION 3.10.  *Let $(Z^\varepsilon)$ be a sequence of continuous and $\mathbb{H}$-adapted processes, and $N$ a continuous $\mathbb{H}$-local martingale. Suppose that $(Z^\varepsilon)$ converges to zero in $\mathcal{A}(\mathbb{H})$. Then for every $h$ in $\mathcal{L}_N^2(\mathbb{H})$, and $Y = \int_0^\cdot h_s \, dN_s$,*

$$\lim_{\varepsilon \to 0} \frac{1}{\varepsilon} \int_0^\cdot Z_s^\varepsilon (Y_{s+\varepsilon} - Y_s) \, ds = 0 \qquad ucp.$$

PROOF.  Since the convergence in probability is equivalent to existence of subsequences convergent to zero almost surely, it is not reductive to suppose that $(Z^\varepsilon)$ converges uniformly to zero, almost surely. We set, for every $k$ in $\mathbb{N}^*$,

$$\Omega_k = \left\{ \omega \in \Omega, \ s.t. \ \sup_{0 \leq s \leq 1} |Z_s^\varepsilon| \leq k, \forall \varepsilon \leq k^{-1} \right\}$$

and

$$Z^{\varepsilon,k} = Z^\varepsilon I_{\{\sup_{0 \leq u \leq \cdot} |Z_u^\varepsilon| \leq k\}}.$$

Then it is sufficient to show that

$$\lim_{\varepsilon \to 0} C_\varepsilon^k = \lim_{\varepsilon \to 0} \frac{1}{\varepsilon} \int_0^\cdot Z_s^{\varepsilon,k} (Y_{s+\varepsilon} - Y_s) \, ds = 0 \qquad ucp, \ \forall k \in \mathbb{N}^*.$$

Let $k$ be fixed. Thanks to adaptedness of the process $Z^{\varepsilon,k}$, we can write

$$(15) \qquad C_\varepsilon^k = \frac{1}{\varepsilon} \int_0^\cdot \left( \int_s^{s+\varepsilon} Z_s^{\varepsilon,k} h_r \, dN_r \right) ds.$$

Let $(\tau^n)_{n \in \mathbb{N}}$ be a sequence of $\mathbb{H}$-stopping times such that $N^{\tau^n}$, the local martingale $N$ stopped at time $\tau^n$, is a square integrable martingale and

$$\sup_{0 \le s \le \tau^n} |h_s| \le n.$$

Stopping integral (15) to time $\tau^n$, let us apply Exercise 5.17, page 165 of [24] to write

$$\frac{1}{\varepsilon} \int_0^{\cdot \wedge \tau^n} \left( \int_s^{(s+\varepsilon) \wedge \tau^n} Z_s^{\varepsilon,k} h_r \, dN_r^{\tau^n} \right) ds = \int_0^{\cdot \wedge \tau^n} \left( \frac{1}{\varepsilon} \int_{(r-\varepsilon)}^r Z_s^{\varepsilon,k} h_r \, ds \right) dN_r^{\tau^n}.$$

By Proposition 2.74 of [18], we are allowed to take the limit for $n \to +\infty$, and write

$$C_\varepsilon^k = \int_0^{\cdot} \left( \frac{1}{\varepsilon} \int_{r-\varepsilon}^r Z_s^{\varepsilon,k} h_r \, ds \right) dN_r \qquad \text{a.s.}$$

Using Doob and Hölder inequalities, we obtain

$$E\left[ \sup_{t \in [0,1]} |C_\varepsilon^k(t)|^2 \right] \le cE\left[ \int_0^1 \left( \frac{1}{\varepsilon} \int_{r-\varepsilon}^r Z_s^{\varepsilon,k} h_r \, ds \right)^2 d[N]_r \right]$$

$$\le cE\left[ \sup_{s \in [0,1]} |Z_s^{\varepsilon,k}|^2 \int_0^1 h_r^2 \, d[N]_r \right],$$

for some positive constant $c$. The Lebesgue dominated convergence theorem permits us to complete the proof. $\square$

COROLLARY 3.11. *Let $R$ be an $\mathbb{H}$-strongly predictable continuous process. Then for every continuous $\mathbb{H}$-local martingale $N$, and $h$ in $\mathcal{L}_N^2(\mathbb{H})$, $[R, Y] = 0$.*

PROOF. It has to be shown that $\left( \frac{1}{\varepsilon} \int_0^{\cdot} Z_s^{\varepsilon}(Y_{s+\varepsilon} - Y_s) \, ds \right)$ converges to zero *ucp*, with $Z^{\varepsilon} = R_{\varepsilon + \cdot} - R$. Since $R$ is $\mathbb{H}$-strongly predictable, $Z^{\varepsilon}$ is definitely $\mathbb{H}$-adapted. Moreover, the continuity of $R$ insures the uniformly convergence to zero, almost surely, of $Z^{\varepsilon}$. Proposition 3.10 leads to the conclusion. $\square$

REMARK 3.12. $[R, Y]$ is zero either for pathwise regularity or for probabilistic reasons. The first situation arises if $R$ has zero quadratic variation when its paths are, for instance, Hölder continuous with parameter $\gamma > \frac{1}{2}$. The second (probabilistic) reason arises, for example, when $R$ is strongly predictable as Corollary 3.11 shows.

We go on defining and discussing some properties of weak Dirichlet processes.

DEFINITION 3.13. An $\mathbb{H}$-*weak Dirichlet* process is the sum of a continuous $\mathbb{H}$-local martingale $M$ and a continuous process $Q$ such that $[Q, N] = 0$, for every $\mathbb{H}$-local martingale $N$.

Corollary 3.11 directly implies the following.

COROLLARY 3.14. *An $\mathbb{H}$-strongly predictable continuous process $R$ is an $\mathbb{H}$-weak Dirichlet process.*

Proposition 3.9 permits us to better specify the nature of such processes with respect to Brownian filtrations, as pointed out in the corollary below.

COROLLARY 3.15. *Suppose that $W$ is a Brownian motion on $(\Omega, \mathcal{F}, P)$. Let $\mathbb{H}$ be its natural filtration augmented by the $P$ null sets. An $\mathbb{H}$-adapted, with finite quadratic variation and continuous process $D$ is an $\mathbb{H}$-Dirichlet process if and only if it is the sum of a continuous $\mathbb{H}$-local martingale $M$ and a finite quadratic variation process $Q$, continuous, $\mathbb{H}$-adapted and such that $[Q, W] = 0$.*

PROOF. Necessity is obvious. Suppose that $D$ is the sum of an $\mathbb{H}$-local martingale $M$ and a continuous process $Q$, with finite quadratic variation, $\mathbb{H}$-adapted and such that $[Q, W] = 0$. Let $N$ be an $\mathbb{H}$-local martingale. Then there exists a process $h$ in $\mathcal{L}_W^2(\mathbb{H})$ such that $N = N_0 + \int_0^\cdot h_s \, dW_s$. By Proposition 3.9, $[Q, N] = \int_0^\cdot h_s \, d[Q, W]_s = 0$. $\square$

THEOREM 3.16. *Let $(X(t, x), 0 \le t \le 1, x \in \mathbb{R})$ be the sum of a $C^1$ $\mathbb{H}$-Itô-martingale field of the form* (10) *and an $\mathbb{H}$-strict zero quadratic variation process $Z$ in $C^{0,1}([0, 1] \times \mathbb{R})$. Let $\xi$ be such that the vector $(\xi, N^1, \ldots, N^n)$ satisfies the hypothesis* (D) *with respect to the filtration $\mathbb{H}$. Then for any semimartingale of the form $Y = \sum_{i=1}^n \int_0^\cdot h_s^i \, dN_s^i$, with $h^i$ in $\mathcal{L}_{N^i}^2(\mathbb{H})$ for every $i = 1, \ldots, n$, it holds:*

$$[X(\cdot, \xi), Y] = \sum_{i=1}^n \int_0^\cdot \partial_x X(s, \xi_s) h_s^i \, d[\xi, N^i]_s$$
$$+ \sum_{i,j=1}^n \int_0^\cdot a^j(s, \xi_s) h_s^i \, d[N^i, N^j]_s.$$

*In particular, $[X(\cdot, \xi), Y]$ has bounded variation.*

REMARK 3.17. In [10] the authors explore the existence of mutual brackets of Itô fields, and so it could appear natural to do the same in this context. However, it is clear that in this case such a bracket cannot exist unless $R$ is a finite quadratic variation process.

PROOF OF THEOREM 3.16. We suppose, for simplicity of notation, that $n = 1$, and we denote with $h$ the process $h^1$. We have to study the convergence $ucp$ of

$$C_\varepsilon(t) = \frac{1}{\varepsilon} \int_0^t (X(s+\varepsilon, \xi_{s+\varepsilon}) - X(s, \xi_s))(Y_{s+\varepsilon} - Y_s)\, ds.$$

We have

$$
\begin{aligned}
C_\varepsilon(t) &= \frac{1}{\varepsilon} \int_0^t (X(s+\varepsilon, \xi_{s+\varepsilon}) - X(s+\varepsilon, Q_s + R_{s+\varepsilon}))(Y_{s+\varepsilon} - Y_s)\, ds \\
&\quad + \frac{1}{\varepsilon} \int_0^t (X(s+\varepsilon, Q_s + R_{s+\varepsilon}) - X(s, Q_s + R_{s+\varepsilon}))(Y_{s+\varepsilon} - Y_s)\, ds \\
&\quad + \frac{1}{\varepsilon} \int_0^t (X(s, Q_s + R_{s+\varepsilon}) - X(s, \xi_s))(Y_{s+\varepsilon} - Y_s)\, ds \\
&= J_\varepsilon^1(t) + J_\varepsilon^2(t) + J_\varepsilon^3(t).
\end{aligned}
$$

For $J_\varepsilon^1(t)$, we use the Taylor type formula

$$
\begin{aligned}
X(s+\varepsilon, \xi_{s+\varepsilon}) &- X(s+\varepsilon, Q_s + R_{s+\varepsilon}) \\
&= \partial_x X(s, \xi_s)(Q_{s+\varepsilon} - Q_s) + \rho(s, \varepsilon)(Q_{s+\varepsilon} - Q_s),
\end{aligned}
$$

with

$$\rho(s, \varepsilon) = \int_0^1 [\partial_x X(s+\varepsilon, \lambda(Q_{s+\varepsilon} - Q_s) + (Q_s + R_{s+\varepsilon})) - \partial_x X(s, \xi_s)]\, d\lambda,$$

to get

$$
\begin{aligned}
J_\varepsilon^1(t) &= \frac{1}{\varepsilon} \int_0^t \partial_x X(s, \xi_s)(Q_{s+\varepsilon} - Q_s)(Y_{s+\varepsilon} - Y_s)\, ds \\
&\quad + \frac{1}{\varepsilon} \int_0^t \rho(s, \varepsilon)(Q_{s+\varepsilon} - Q_s)(Y_{s+\varepsilon} - Y_s)\, ds.
\end{aligned}
$$

Since $h$ is continuous and $\mathbb{H}$-adapted, it is progressively measurable and almost surely bounded. By Proposition 3.9, $(Q, \int_0^{\cdot} h_s\, dN_s)$ has all its mutual brackets, and so by Remark 2.5.2, the first term converges $ucp$ to

$$\int_0^{\cdot} \partial_x X(s, \xi_s) h_s\, d[Q, N]_s,$$

while the second term has limit equal to zero $ucp$ since both $Q$ and $Y$ have *finite quadratic variation*.

We consider the term $J^2(t)$. Thanks to the hypothesis $(\mathcal{D})$, the process

$$(Q_s + R_{s+\varepsilon}, 0 \le s \le 1)$$

is $\mathbb{H}$-adapted for every $\varepsilon \leq \delta$. Then we can write, for every $\varepsilon \leq \delta$,

$$
\begin{aligned}
J_\varepsilon^2(t) = {} & \frac{1}{\varepsilon} \int_0^t \left( \int_s^{s+\varepsilon} (a(r, Q_s + R_{s+\varepsilon}) - a(r, \xi_r)) \, dN_r \right) (Y_{s+\varepsilon} - Y_s) \, ds \\
& + \frac{1}{\varepsilon} \int_0^t \left( \int_s^{s+\varepsilon} a(r, \xi_r) \, dN_r \right) (Y_{s+\varepsilon} - Y_s) \, ds \\
& + \frac{1}{\varepsilon} \int_0^t (Z(s+\varepsilon, Q_s + R_{s+\varepsilon}) - Z(s, Q_s + R_{s+\varepsilon}))(Y_{s+\varepsilon} - Y_s) \, ds.
\end{aligned}
$$

The second term converges *ucp* by definition to

$$
\left[ \int_0^\cdot a(s, \xi_s) \, dN_s, Y \right] = \int_0^t h_s a(s, \xi_s) \, d[N, N]_s,
$$

while, using the Hölder inequality and the fact that $Z$ is a *strict zero quadratic variation* process, it is possible to show that the last term converges to zero *ucp*. Again by the Hölder inequality, the first term converges to zero *ucp* if

$$
\lim_{\varepsilon \to 0} \frac{1}{\varepsilon} \int_0^1 \left( \int_s^{s+\varepsilon} (a(r, Q_s + R_{s+\varepsilon}) - a(r, \xi_r)) \, dN_r \right)^2 ds = 0, \qquad \text{in probability.}
$$

This can be proved with techniques already used for the convergence to zero of the term $I_\varepsilon^2$ in the proof of Proposition 3.4. Regarding the term $J^3$, we apply Proposition 3.10 to the sequence of processes $(X(\cdot, Q + R_{\cdot+\varepsilon}) - X(\cdot, \xi))$, the local martingale $N$ and the process $h$, which let us conclude that $J^3$ converges to zero *ucp*. $\quad \square$

Using similar arguments to those of the previous proposition, one can prove the following.

PROPOSITION 3.18. *Let $\beta$ be in $C^{0,1}([0,1] \times \mathbb{R})$ and $(\xi, N^1, \ldots, N^n)$ be a vector of continuous processes satisfying the hypothesis (D) with respect to $\mathbb{H}$. Then for every semimartingale of the form*

$$
Y = \sum_{i=1}^n \int_0^\cdot h_s^i \, dN_s^i,
$$

*with $h^i$ in $\mathcal{L}_{N^i}^2(\mathbb{H})$ for every $i = 1, \ldots, n$, $[\beta(\cdot, \xi), Y]$ exists and*

$$
(16) \qquad [\beta(\cdot, \xi), Y] = \sum_{i=1}^n \int_0^\cdot h_s^i \partial_x \beta(s, \xi_s) \, d[\xi, N^i]_s.
$$

*In particular, $[\beta(\cdot, \xi), Y]$ has bounded variation.*

COROLLARY 3.19. *Let $(\xi, N^1, \ldots, N^n)$ be a vector of continuous processes satisfying the hypothesis ($\mathcal{D}$) with respect to $\mathbb{H}$. Let $X = (X(t,x), 0 \leq t \leq 1, x \in \mathbb{R})$ and $Z = (Z(t,x), 0 \leq t \leq 1, x \in \mathbb{R})$ be either functions in $C^{0,1}([0,1] \times \mathbb{R})$ or $C^1$ $\mathbb{H}$-Itô-semimartingale fields of the form* (13). *Then for every semimartingale of the form*

$$Y = \sum_{i=1}^{n} \int_0^{\cdot} h_s^i \, dN_s^i,$$

*with $h^i$ in $\mathcal{L}^2_{N^i}(\mathbb{H})$ for every $i = 1, \ldots, n$, it holds*

$$\int_0^{\cdot} X(s, \xi_s) \, d^\circ \left( \int_0^s Z(r, \xi_r) \, d^\circ Y_r \right) = \int_0^{\cdot} (XZ)(s, \xi_s) \, d^\circ Y_s.$$

PROOF. The corollary is a consequence of Proposition 3.16 and the decomposition of the symmetric integral into a classical stochastic integral plus a half covariation as specified in Remark 2.7.2. $\square$

3.3. *Itô–Wentzell formula.* The process $\xi$ in this subsection will be always supposed to be a strong cubic variation process.

PROPOSITION 3.20. *Let $(X(t,x), 0 \leq t \leq 1, x \in \mathbb{R})$ be a $C^3$ $\mathbb{H}$-Itô-semimartingale field of the form* (13). *Let $(\xi, N^1, \ldots, N^n)$ be a vector of continuous processes satisfying the hypothesis ($\mathcal{D}$) with respect to $\mathbb{H}$. Then the symmetric integral $\int_0^{\cdot} \partial_x X(s, \xi_s) \, d^\circ \xi_s$ exists and*

$$X(\cdot, \xi) = X(0, \xi_0) + \sum_{i=1}^{n} \int_0^{\cdot} a^i(s, \xi_s) \, dN_s^i + \sum_{j=1}^{m} \int_0^{\cdot} b^j(s, \xi_s) \, dV_s^j$$

$$+ \int_0^{\cdot} \partial_x X(s, \xi_s) \, d^\circ \xi_s + \frac{1}{2} \sum_{i=1}^{n} \int_0^{\cdot} \partial_x a^i(s, \xi_s) \, d[N^i, \xi]_s$$

$$- \frac{1}{12} \int_0^{\cdot} \partial_x^{(3)} X(s, \xi_s) \, d[\xi, \xi, \xi]_s.$$

PROOF. We suppose $n = m = 1$, and we make the usual simplification in the notation of the Itô field considered. By continuity of the process $X(\cdot, \xi)$, the sequence of processes

$$\frac{1}{\varepsilon} \int_0^t (X(s + \varepsilon, \xi_{s+\varepsilon}) - X(s, \xi_s)) \, ds$$

converges almost surely to $(X(t, \xi_t) - X(0, \xi_0))$. In particular,

$$X(t, \xi_t) - X(0, \xi_0) = \lim_{\varepsilon \to 0} \frac{1}{\varepsilon} \int_0^t (X(s + \varepsilon, \xi_{s+\varepsilon}) - X(s + \varepsilon, \xi_s)) \, ds$$

$$+ \lim_{\varepsilon \to 0} \frac{1}{\varepsilon} \int_0^t (X(s+\varepsilon, \xi_s) - X(s, \xi_s)) \, ds$$

$$= \lim_{\varepsilon \to 0} I_\varepsilon^1(t) + \lim_{\varepsilon \to 0} I_\varepsilon^2(t),$$

if the two limits on the right-hand side of the previous equality exist. Applying substitution arguments and interchanging the integrals with respect to time, the semimartingales $N$ and $V$, $I_\varepsilon^2(t)$ converges *ucp* to

$$\int_0^\cdot a(s, \xi_s) \, dN_s + \int_0^\cdot b(s, \xi_s) \, dV_s.$$

Since $X(\cdot, x)$ is differentiable till order three with respect to $x$, we can write

$$(17) \quad \begin{aligned} &X(s+\varepsilon, \xi_{s+\varepsilon}) \\ &= X(s+\varepsilon, \xi_s) + \partial_x X(s+\varepsilon, \xi_s)(\xi_{s+\varepsilon} - \xi_s) \\ &\quad + \tfrac{1}{2}\partial_x^{(2)} X(s+\varepsilon, \xi_s)(\xi_{s+\varepsilon} - \xi_s)^2 \\ &\quad + \tfrac{1}{6}\partial_x^{(3)} X(s+\varepsilon, \xi_s)(\xi_{s+\varepsilon} - \xi_s)^3 + \rho(\xi_s, \xi_{s+\varepsilon})(\xi_{s+\varepsilon} - \xi_s)^3 \end{aligned}$$

and

$$(18) \quad \begin{aligned} &X(s+\varepsilon, \xi_s) \\ &= X(s+\varepsilon, \xi_{s+\varepsilon}) + \partial_x X(s+\varepsilon, \xi_{s+\varepsilon})(\xi_s - \xi_{s+\varepsilon}) \\ &\quad + \tfrac{1}{2}\partial_x^{(2)} X(s+\varepsilon, \xi_{s+\varepsilon})(\xi_s - \xi_{s+\varepsilon})^2 \\ &\quad + \tfrac{1}{6}\partial_x^{(3)} X(s+\varepsilon, \xi_{s+\varepsilon})(\xi_s - \xi_{s+\varepsilon})^3 + \rho(\xi_{s+\varepsilon}, \xi_s)(\xi_s - \xi_{s+\varepsilon})^3, \end{aligned}$$

with $\lim_{\varepsilon \to 0} \rho(\xi_s, \xi_{s+\varepsilon}) = \lim_{\varepsilon \to 0} \rho(\xi_{s+\varepsilon}, \xi_s) = 0$, almost surely. By subtracting these two quantities and integrating over $[0, t]$, we get

$$\begin{aligned} I_\varepsilon^1(t) &= \frac{1}{2\varepsilon} \int_0^t (\partial_x X(s+\varepsilon, \xi_s) + \partial_x X(s+\varepsilon, \xi_{s+\varepsilon}))(\xi_{s+\varepsilon} - \xi_s) \, ds \\ &\quad - \frac{1}{4\varepsilon} \int_0^t (\partial_x^{(2)} X(s+\varepsilon, \xi_{s+\varepsilon}) - \partial_x^{(2)} X(s+\varepsilon, \xi_s))(\xi_{s+\varepsilon} - \xi_s)^2 \, ds \\ &\quad + \frac{1}{12\varepsilon} \int_0^t (\partial_x^{(3)} X(s+\varepsilon, \xi_s) + \partial_x^{(3)} X(s+\varepsilon, \xi_{s+\varepsilon}))(\xi_{s+\varepsilon} - \xi_s)^3 \, ds \\ &\quad + \frac{1}{2\varepsilon} \int_0^t (\rho(\xi_s, \xi_{s+\varepsilon}) + \rho(\xi_{s+\varepsilon}, \xi_s))(\xi_{s+\varepsilon} - \xi_s)^3 \, ds \\ &= J_\varepsilon^1(t) + J_\varepsilon^2(t) + J_\varepsilon^3(t) + J_\varepsilon^4(t). \end{aligned}$$

Since $\xi$ is a *strong cubic variation* process, $J_\varepsilon^4$ converges to zero *ucp*. $J_\varepsilon^2$ converges *ucp* to

$$-\tfrac{1}{4}[\partial_x^{(2)} X(\cdot, \xi), \xi, \xi].$$

In fact,

$$J_\varepsilon^2(t) = -\frac{1}{4\varepsilon} \int_0^t (\partial_x^{(2)} X(s+\varepsilon, \xi_{s+\varepsilon}) - \partial_x^{(2)} X(s, \xi_s))(\xi_{s+\varepsilon} - \xi_s)^2 \, ds$$

$$+ \frac{1}{4\varepsilon} \int_0^t (\partial_x^{(2)} X(s+\varepsilon, \xi_s) - \partial_x^{(2)} X(s, \xi_s))(\xi_{s+\varepsilon} - \xi_s)^2 \, ds.$$

The first term converges $ucp$ to

$$-\tfrac{1}{4}[\partial_x^{(2)} X(\cdot, \xi), \xi, \xi] = -\tfrac{1}{4} \int_0^\cdot \partial_x^{(3)} X(s, \xi_s) \, d[\xi, \xi, \xi]_s,$$

since $\partial_x^2 X(\cdot, x)$ is a $C^1$ $\mathbb{H}$-Itô-semimartingale field and Proposition 3.4 can be applied. The second term converges to zero $ucp$. In fact, by the Hölder inequality, its absolute value is bounded by

$$\tfrac{1}{4} \left( \frac{1}{\varepsilon} \int_0^1 |\partial_x^{(2)} X(s+\varepsilon, \xi_s) - \partial_x^{(2)} X(s, \xi_s)|^3 \, ds \right)^{1/3} ds \|[\xi, \xi, \xi]\|_\varepsilon^{2/3}.$$

Since $\partial_x^{(2)} X$ is a $C^1$ Itô-semimartingale field, the first factor of the product can be shown to converge to zero in probability, using tools already developed in the proof of Proposition 3.4 for the term $\int_0^\cdot |B(s, \varepsilon)|^3 \, ds$. The term $J_\varepsilon^3$ can be written as

$$\frac{1}{12\varepsilon} \int_0^t (\partial_x^{(3)} X(s+\varepsilon, \xi_{s+\varepsilon}) + \partial_x^{(3)} X(s+\varepsilon, \xi_s) - 2\partial_x^{(3)} X(s, \xi_s))(\xi_{s+\varepsilon} - \xi_s)^3 \, ds$$

$$+ \frac{1}{6\varepsilon} \int_0^t \partial_x^{(3)} X(s, \xi_s)(\xi_{s+\varepsilon} - \xi_s)^3 \, ds.$$

By Remark 2.5.2, the second term converges $ucp$ to $\frac{1}{6} \int_0^\cdot \partial_x^{(3)} X(s, \xi_s) \, d[\xi, \xi, \xi]_s$, while the first term converges to zero 0 a.s., since $\xi$ is a finite strong cubic variation, and both $\partial_x^{(3)} X$ and $\xi$ are continuous. Finally,

$$J_\varepsilon^1 = \frac{1}{2\varepsilon} \int_0^t (\partial_x X(s, \xi_s) + \partial_x X(s+\varepsilon, \xi_{s+\varepsilon}))(\xi_{s+\varepsilon} - \xi_s) \, ds$$

$$+ \frac{1}{2\varepsilon} \int_0^t (\partial_x X(s+\varepsilon, \xi_s) - \partial_x X(s, \xi_s))(\xi_{s+\varepsilon} - \xi_s) \, ds.$$

The second term can be decomposed in the following way:

$$\frac{1}{2\varepsilon} \int_0^t (\partial_x X(s+\varepsilon, \xi_s) - \partial_x X(s, \xi_s))(\xi_{s+\varepsilon} - \xi_s) \, ds$$

$$= \frac{1}{2\varepsilon} \int_0^t \left( \int_s^{s+\varepsilon} (\partial_x a(r, \xi_s) - \partial_x a(r, \xi_r)) \, dN_r \right)(Q_{s+\varepsilon} - Q_s) \, ds$$

$$+ \frac{1}{2\varepsilon} \int_0^t \left( \int_s^{s+\varepsilon} (\partial_x a(r, \xi_s) - \partial_x a(r, \xi_r))(R_{s+\varepsilon} - R_s) \, dN_r \right) ds$$

$$+ \frac{1}{2\varepsilon} \int_0^t \left( \int_s^{s+\varepsilon} \partial_x a(r, \xi_r) \, dN_r \right) (\xi_{s+\varepsilon} - \xi_s) \, ds$$

$$+ \frac{1}{2\varepsilon} \int_0^t (Z(s+\varepsilon, \xi_s) - Z(s, \xi_s))(\xi_{s+\varepsilon} - \xi_s) \, ds,$$

with $Z = \int_0^\cdot \partial_x b(s, \cdot) \, dV_s$. The first term of the sum converges to zero *ucp* by the Hölder inequality, since $Q$ is a finite quadratic variation process and

$$\lim_{\varepsilon \to 0} \int_0^1 \frac{1}{\varepsilon} \left( \int_s^{s+\varepsilon} (\partial_x a(r, \xi_s) - \partial_x a(r, \xi_r)) \, dN_r \right)^2 ds = 0, \qquad \text{in probability.}$$

By Proposition 3.16,

$$\lim_{\varepsilon \to 0} \frac{1}{2\varepsilon} \int_0^t \left( \int_s^{s+\varepsilon} \partial_x a(r, \xi_r) \, dN_r \right) (\xi_{s+\varepsilon} - \xi_s) \, ds$$

$$= \frac{1}{2} \int_0^\cdot \partial_x a(s, \xi_s) \, d[N, \xi]_s, \qquad ucp.$$

The second term can be shown to converge to zero by arguments used in the proof of Proposition 3.10, while the last term converges to zero *ucp* since $Z$ is an $\mathbb{H}$-strict zero $p$-variation process, for every $p > 1$. As a consequence of this, the first term of $J_\varepsilon^1$ is forced to converge to

$$\int_0^\cdot \partial_x X(s, \xi_s) \, d^\circ \xi_s,$$

and we get the result. $\square$

3.4. *Existence of symmetric integrals and chain-rule formulae.* The process $\xi$ in this subsection will be always supposed to be a strong cubic variation process.

DEFINITION 3.21. We will denote with $\mathcal{C}_\xi^k(\mathbb{H})$ the set of all processes of the form

$$Z_t = X(t, \xi_t),$$

being $X$ a $C^k$ $\mathbb{H}$-Itô-semimartingale field driven by the vector of local martingales $(N^1, \ldots, N^n)$, such that the vector $(\xi, N^1, \ldots, N^n)$ satisfies the hypothesis ($\mathcal{D}$), with respect to the filtration $\mathbb{H}$.

REMARK 3.22. The set $\mathcal{C}_\xi^k(\mathbb{H})$ is an algebra (apply classical Itô formula).

REMARK 3.23. 1. A process $Z$ belongs to $\mathcal{C}_\xi^3(\mathbb{H})$ if and only if there exist an $\mathcal{H}_0$-measurable random variable $Z_0$, a vector of $\mathbb{H}$-adapted processes $(N^1, \ldots, N^n)$ such that $(\xi, N^1, \ldots, N^n)$ satisfies the hypothesis $(\mathcal{D})$ with respect to $\mathbb{H}$, a vector of $\mathbb{H}$-adapted stochastic processes $(h^1, \ldots, h^n)$, and a process $\gamma$ in $\mathcal{C}^2_\xi(\mathbb{H})$, such that

$$Z = Z_0 + \int_\circ^\cdot \gamma_s \, d^\circ \xi_s + \sum_{i=1}^n \int_0^\cdot h_s \, dN^i_s.$$

The statement is a direct consequence of the Itô–Wentzell formula.

2. Combining Remark 2.5, the *reversed* Itô–Wentzell formula, and Proposition 3.4, it is possible to prove that if $\gamma^1, \gamma^2$ and $\gamma^3$ belong to $\mathcal{C}^2_\xi(\mathbb{H})$, then

$$\left[ \int_0^\cdot \gamma^1_s d^\circ \xi_s, \int_0^\cdot \gamma^2_s \, d^\circ \xi_s, \int_0^\cdot \gamma^3_s \, d^\circ \xi_s \right] = \int_0^\cdot \gamma^1_s \gamma^2_s \gamma^3_s [\xi, \xi, \xi]_s.$$

3. A significant example of the class $\mathcal{C}^3_\xi(\mathbb{H})$ is given by the following. Let $W = (W^1, \ldots, W^n)$ be a $n$-dimensional Brownian motion on $(\Omega, \mathcal{F}, P)$ with respect to its natural filtration $\mathbb{H}$ augmented by the $P$ null sets. Let $\xi$ be a strong cubic variation process such that the vector $(\xi, W^1, \ldots, W^n)$ satisfies the hypothesis $(\mathcal{D})$ with respect to $\mathbb{H}$. Then the set $\mathcal{C}^3_\xi(\mathbb{H})$ coincides with the processes of the form

$$Z = Z_0 + \int_0^\cdot \gamma_s \, d^\circ \xi_s + L,$$

where $\gamma$ is in $\mathcal{C}^2_\xi(\mathbb{H})$ and $L$ is an $\mathbb{H}$-semimartingale. This holds since every $\mathbb{H}$-local martingale, zero at $t = 0$, admits a representation as a stochastic integral with respect to $W$.

PROPOSITION 3.24. *For every $Z$ in $\mathcal{C}^2_\xi(\mathbb{H})$ and $U$ in $\mathcal{C}^3_\xi(\mathbb{H})$, the symmetric integral*

$$\int_0^\cdot Z_s \, d^\circ U_s,$$

*exists and belongs to $\mathcal{C}^2_\xi(\mathbb{H})$. If $Z_t = Y(t, \xi_t)$ and $U_t = X(t, \xi_t)$, where $X(\cdot, x)$ and $Y(\cdot, x)$ have representations*

$$(19) \qquad X(\cdot, x) = X_0(x) + \sum_{i=1}^n \int_0^\cdot a^i(s, x) \, dN^i_s + \sum_{j=1}^m \int_0^\cdot b^j(s, x) \, dV^j_s$$

*and*

$$(20) \qquad Y(\cdot, x) = Y_0(x) + \sum_{i=1}^n \int_0^\cdot \bar{a}^i(s, x) \, dN^i_s + \sum_{j=1}^m \int_0^\cdot \bar{b}^j(s, x) \, dV^j_s,$$

*then*

$$\int_0^{\cdot} Z_s \, d^{\circ} U_s = \sum_{i=1}^{n} \int_0^{\cdot} (Y a^i)(s, \xi_s) \, dN_s^i + \sum_{j=1}^{m} \int_0^{\cdot} (Y b^j)(s, \xi_s) \, dV_s^j$$

$$+ \int_0^{\cdot} (Y \partial_x X)(s, \xi_s) \, d^{\circ} \xi_s + \frac{1}{2} \sum_{i=1}^{n} \int_0^{\cdot} \partial_x (Y a^i)(s, \xi_s) \, d[N^i, \xi]_s$$

$$+ \frac{1}{2} \sum_{i,j=1}^{n} \int_0^{\cdot} (a^j \bar{a}^i)(s, \xi_s) \, d[N^i, N^j]_s$$

$$- \frac{1}{12} \int_0^{\cdot} ((3 \partial_x^{(2)} X)(\partial_x Y) + (\partial_x^{(3)} X) Y)(s, \xi_s) \, d[\xi, \xi, \xi]_s.$$

PROOF. We restrict ourselves to the case $n = m = 1$, and we denote $a^1 = a, \bar{a}^1 = \bar{a}$. We have to investigate the convergence of

$$C_{\varepsilon}(t) = \frac{1}{2\varepsilon} \int_0^t (Z_{s+\varepsilon} + Z_s)(U_{s+\varepsilon} - U_s) \, ds$$

$$= \frac{1}{2\varepsilon} \int_0^t (Z_{s+\varepsilon} + Z_s)(X(s+\varepsilon, \xi_{s+\varepsilon}) - X(s+\varepsilon, \xi_s)) \, ds$$

$$+ \frac{1}{2\varepsilon} \int_0^t (Z_{s+\varepsilon} + Z_s)(X(s+\varepsilon, \xi_s) - X(s, \xi_s)) \, ds$$

$$= I_{\varepsilon}^1(t) + I_{\varepsilon}^2(t).$$

As concerns the second term, we can write

$$I_{\varepsilon}^2(t) = \frac{1}{2\varepsilon} \int_0^t (Z_{s+\varepsilon} - Z_s)(X(s+\varepsilon, \xi_s) - X(s, \xi_s)) \, ds$$

$$+ \frac{1}{\varepsilon} \int_0^t Z_s (X(s+\varepsilon, \xi_s) - X(s, \xi_s)) \, ds.$$

Using techniques already introduced in the previous section and in Proposition 3.16, one can show that these two terms converge, respectively, *ucp* to

$$\frac{1}{2} \left[ Y(\cdot, \xi), \int_0^{\cdot} a(r, \xi_r) \, dN_r + \int_0^{\cdot} b(r, \xi_r) \, dV_r \right]$$

$$= \frac{1}{2} \int_0^{\cdot} ((\partial_x Y) a)(s, \xi_s) \, d[N, \xi]_s + \frac{1}{2} \int_0^{\cdot} (\bar{a} a)(s, \xi_s) \, d[N, N]_s$$

and $\int_0^{\cdot} Z_s a(s, \xi_s) \, dN_s + \int_0^{\cdot} Z_s b(s, \xi_s) \, dV_s$.

We consider the first term. To this extent, for every $s$ in $[0, 1]$, we multiply equalities (17) and (18), respectively, by $Z_s$ and $Z_{s+\varepsilon}$ to get

$$I_{\varepsilon}^1(t) = \frac{1}{2\varepsilon} \int_0^t (\partial_x X(s+\varepsilon, \xi_s) Z_s + \partial_x X(s+\varepsilon, \xi_{s+\varepsilon}) Z_{s+\varepsilon})(\xi_{s+\varepsilon} - \xi_s) \, ds$$

$$- \frac{1}{4\varepsilon} \int_0^t (\partial_x^{(2)} X(s+\varepsilon, \xi_{s+\varepsilon}) Z_{s+\varepsilon} - \partial_x^{(2)} X(s+\varepsilon, \xi_s) Z_s)(\xi_{s+\varepsilon} - \xi_s)^2 \, ds$$

$$+ \frac{1}{12\varepsilon} \int_0^t (\partial_x^{(3)} X(s+\varepsilon, \xi_s) Z_s + \partial_x^{(3)} X(s+\varepsilon, \xi_{s+\varepsilon}) Z_{s+\varepsilon})(\xi_{s+\varepsilon} - \xi_s)^3 \, ds$$

$$+ \frac{1}{2\varepsilon} \int_0^t (\rho(\xi_s, \xi_{s+\varepsilon}) Z_s + \rho(\xi_{s+\varepsilon} Z_{s+\varepsilon}, \xi_s) Z_s)(\xi_{s+\varepsilon} - \xi_s)^3 \, ds.$$

The proof follows the same outlines of the calculus already performed in the proof of the Itô–Wentzell formula for the term $I_\varepsilon^1(t)$. The Itô–Wentzell formula is indeed a particular case of this result ($Z = 1$). The only difference, here, is that the symmetric integral $\int_0^\cdot \partial_x X(s, \xi_s) Z_s \, d^\circ \xi_s$ exists since $\partial_x X(\cdot, x) Z$ is still a $C^2$ $\mathbb{H}$-Itô-semimartingale field, and for such a field, the existence was already proved before. Then, similarly, we obtain

$$\lim_{\varepsilon \to 0} I_\varepsilon^1(t) = \int_0^t Z_s \partial_x X(s, \xi_s) \, d^\circ \xi_s + \frac{1}{2} \int_0^t Z_s \partial_x a(s, \xi_s) \, d[N, \xi]_s$$

$$- \frac{1}{4} [\partial_x^{(2)} X(\cdot, \xi) Z, \xi, \xi]_t + \frac{1}{6} \int_0^t Z_s \partial_x^{(3)} X(s, \xi_s) \, d[\xi, \xi, \xi]_s \qquad ucp.$$

The conclusion follows applying Proposition 3.4 to get the equality

$$[\partial_x^{(2)} X(\cdot, \xi) Z, \xi, \xi]_t = \int_0^t (Z_s \partial_x^{(3)} X(s, \xi_s) + \partial_x^{(2)} X \partial_x Y(s, \xi_s)) \, d[\xi, \xi, \xi]_s,$$

which leads to the result. $\square$

PROPOSITION 3.25. *Let $Z$ and $U$ be in $\mathcal{C}_\xi^2(\mathbb{H})$, with $Z_t = Y(t, \xi_t)$ and $U_t = X(t, \xi_t)$, where $X(\cdot, x)$ and $Y(\cdot, x)$ have representations* (19) *and* (20). *Then the symmetric integral*

$$\int_0^\cdot Z_s \, d^\circ \left( \int_0^s U(r) \, d^\circ \xi_r \right)$$

*exists and*

$$\int_0^\cdot Z_s \, d^\circ \left( \int_0^s U_r \, d^\circ \xi_r \right) = \int_0^\cdot Z_s U_s \, d^\circ \xi_s - \frac{1}{4} \int_0^\cdot ((\partial_x X)(\partial_x Y))(s, \xi_s) \, d[\xi, \xi, \xi]_s.$$

PROOF. We consider the field $(X^*(t, x), 0 \le t \le 1, x \in \mathbb{R})$ so defined

$$X^*(t, x) = \int_0^x X(t, z) \, dz.$$

Clearly, the $X^*$ is a $C^3$ $\mathbb{H}$-Itô-semimartingale field, so the Itô–Wentzell formula can be applied to write

$$\int_0^t X(s, \xi_s) \, d^\circ \xi_s$$

$$= X^*(t, \xi_t) - \sum_{i=1}^{n} \int_0^t a^{i,*}(s, \xi_s) \, dN_s^i - \sum_{j=1}^{m} \int_0^t b^{j,*}(s, \xi_s) \, dV_s^j$$

$$- \frac{1}{2} \sum_{i=1}^{n} \int_0^t a^i(s, \xi_s) \, d[\xi, N^i]_s + \frac{1}{12} \int_0^t \partial_x^{(2)} X(s, \xi_s) \, d[\xi, \xi, \xi]_s,$$

where

$$a^{i,*}(t, x) = \int_0^x a^i(t, z) \, dz, \qquad b^{j,*}(t, x) = \int_0^x b^j(t, z) \, dz,$$

for $i = 1, \ldots, n$ and $j = 1, \ldots, m$, are the coefficients comparing in the representation of $X^*$. Since $Y(\cdot, \xi)$ and $X^*(\cdot, \xi)$ are in $\mathcal{C}_\xi^2(\mathbb{H})$ and $\mathcal{C}_\xi^3(\mathbb{H})$, respectively, we can use Propositions 3.16 and 3.24 to conclude. $\square$

## 4. On an SDE driven by a strong cubic variation process and semimartingales.

4.1. *The equation.* On a filtered probability space $(\Omega, \mathcal{F}, \mathbb{F}, P)$, with $\mathbb{F} = (\mathcal{F}_t)_{t \in [0,1]}$, $\mathcal{F}_1 = \mathcal{F}$, let $\xi$, $M$ and $V$ be adapted and, respectively, a *strong cubic variation* process, a local martingale and a bounded variation process. We suppose $\xi_0 = 0$. Let $\sigma, \beta : [0,1] \times \mathbb{R} \to \mathbb{R}$ be continuous functions, $\alpha : [0,1] \times \mathbb{R} \times \Omega \to \mathbb{R}$ be progressively measurable and locally bounded in $x$, uniformly in $t$, almost surely, and $\eta$ be a random variable $\mathcal{F}_0$-measurable.

DEFINITION 4.1. A continuous process $X : \Omega \times [0,1] \to \mathbb{R}$ is called a *solution* to equation

$$(21) \qquad \begin{cases} d^\circ X_t = \sigma(t, X_t)[d^\circ \xi_t + \beta(t, X_t) \, d^\circ M_t + \alpha(t, X_t) \, dV_t], & 0 \leq t \leq 1, \\ X_0 = \eta \end{cases}$$

on $(\Omega, \mathcal{F}, P)$, if:

1. $X_0 = \eta$;
2. $X$ is a *strong cubic variation* process;
3. $[\beta(\cdot, X), M]$ exists and it has bounded variation;
4. for every $\psi$ in $C^{1,\infty}([0,1] \times \mathbb{R})$,

$$\int_0^\cdot \psi(t, X_t) \, d^\circ X_t = \int_0^\cdot (\psi \sigma)(t, X_t)[d^\circ \xi_t + \beta(t, X_t) \, d^\circ M_t + \alpha(t, X_t) \, dV_t]$$

$$- \frac{1}{4} \int_0^\cdot (\partial_x \sigma)(\sigma^2)(\partial_x \psi)(t, X_t) \, d[\xi, \xi, \xi]_t, \qquad \text{a.s.}$$

REMARK 4.2. 1. A solution to equation (21) is a solution to the integral equation

$$X_t = \eta + \int_0^t \sigma(s, X_s) \, d^\circ \xi_s + \int_0^t (\sigma \beta)(s, X_s) \, d^\circ M_s$$

(22)
$$+ \int_0^t (\sigma\alpha)(s, X_s)\, dV_s$$

(consider the case $\psi = 1$).

2. If $X$ is a solution, then property 4 is satisfied for every $\psi$ in $C^{1,2}$ (see [7], Remark 4.2, page 286).

4.2. *Hypotheses on the coefficients.* The construction here used to prove some results about uniqueness and existence of equation (21) is based on the following assumption:

(H$_1$)  $\{(t,x) \in [0,1] \times \mathbb{R}, \text{ s.t. } \sigma(t,x) \neq 0\} = [0,1] \times S = \bigcup_{n=0}^{\infty} ([0,1] \times S^n),$

where $S$ is an open set in $\mathbb{R}$, and thus the countable union of its connected components

$$(S^n = (a_n, b_n), -\infty \leq a_n < b_n \leq +\infty)_{n \in \mathbb{N}}.$$

For every $n$ in $\mathbb{N}$, we define the function $H^n \colon [0,1] \times S^n$:

$$H^n(t,x) = \int_{c_n}^{x} \frac{1}{\sigma(t,z)}\, dz$$

being $c_n$ in $S^n$, and we denote $H(t,x) = \sum_{n=0}^{+\infty} H^n(t,x) I_{[0,1] \times S^n}(t,x)$, for $(t,x)$ in $[0,1] \times S$. We will also need to assume that, for every $t$ in $[0,1]$ and $n$ in $\mathbb{N}$,

(H$_2$)

$$\lim_{(s,x) \to (t,a_n)} \int_x^{c_n} \frac{1}{|\sigma(s,z)|}\, dz = \lim_{(s,x) \to (t,a_n)} |H^n(s,x)| = +\infty,$$

$$\lim_{(s,x) \to (t,b_n)} \int_{c_n}^{x} \frac{1}{|\sigma(s,z)|}\, dz = \lim_{(s,x) \to (t,b_n)} |H^n(s,x)| = +\infty.$$

REMARK 4.3.   1. Assumption (H$_1$) is always verified if $\sigma$ is *autonomous*, that is, if $\sigma(t,x) = \sigma(x)$, for every $0 \leq t \leq 1$.

2. Suppose that $\sigma$ is locally Lipschitz in space, then assumption (H$_2$) is satisfied, for every $n$ in $\mathbb{N}$ such that $-\infty < a_n < b_n < \infty$. In fact, since $\sigma(t, a_n) = \sigma(t, b_n) = 0$ for every $t$, there will be a constant $c > 0$, such that

$$|H^n(t,a)| \geq c(\log(c_n - a_n) - \log(a - a_n)) \qquad \forall\, a \in (a_n, c_n),$$

$$|H^n(t,b)| \geq c(\log(b_n - c_n) - \log(b_n - b)) \qquad \forall\, b \in (c_n, b_n).$$

If $\sigma$ is locally Lipschitz in space, assumption (H$_2$) reduces to verify the nonintegrability condition above only when $a_n$ or $b_n$ are infinity. Even in that case, (H$_2$) is just there to avoid technicalities related to the possible explosion of the solution. As far as uniqueness is concerned, it is not needed.

Under assumption (H$_2$), for every $n$ in $\mathbb{N}$ and $t$ in $[0, 1]$, $H^n(t, \cdot): S^n \to \mathbb{R}$ admits an inverse $K^n(t, \cdot): \mathbb{R} \to S^n$. If $\sigma$ never vanishes, then we will simply denote $K^n$ with $K$. Clearly, for every $n$, $K^n$ is the solution of the following equation:

$$\begin{cases} \partial_y K^n(t, y) = \sigma(t, K^n(t, y)), & (t, y) \in [0, 1] \times \mathbb{R}, \\ K^n(t, 0) = c_n. \end{cases}$$

For every $g: [0, 1] \times S \to \mathbb{R}$, we will denote

$$\widetilde{g}(t, y, \omega) = \sum_{n=0}^{+\infty} I_{\{\eta \in S^n\}}(\omega) g(t, K^n(t, y)), \qquad (t, y, \omega) \in [0, 1] \times \mathbb{R} \times \Omega.$$

4.3. *Some properties on the trajectories of a solution.* The key point of our construction relies on the following property about trajectories of solutions holding if $\sigma$ never vanishes. As we will see, in this case, a solution to equation (21) can be represented in terms of the primitive of $\sigma^{-1}$ which can be defined on $\mathbb{R}$ at every instant. When this is not the case this property will be still true only *locally*, the local character depending on the initial condition $\eta$, and for its description we will need to consider the primitives of $\sigma^{-1}$ on each connected component of $S$.

LEMMA 4.4. *Let $\sigma$ be in $C^{1,2}$, never vanishing and satisfying (H$_2$), $\beta$ be in $C^{0,1}$. Suppose that $X$ is a solution to equation (21) adapted to $\mathbb{F}$. Then*

$$H(\cdot, X) = \xi + N,$$

*where $N$ is the $\mathbb{F}$-semimartingale*

$$N = H(0, \eta) + \int_0^{\cdot} \beta(s, X_s) \, dM_s + \int_0^{\cdot} \alpha(s, X_s) \, dV_s + \int_0^{\cdot} \partial_s H(s, X_s) \, ds$$

$$+ \tfrac{1}{2}[\beta(\cdot, X), M] + \tfrac{1}{12} \int_0^{\cdot} (\sigma \partial_x^{(2)} \sigma + (\partial_x \sigma)^2)(s, X_s) \, d[\xi, \xi, \xi]_s.$$

*Furthermore, if $\sigma$ is autonomous, then the result still holds even if $X$ fulfills property 4 of Definition 4.1, only for autonomous functions $\psi$.*

PROOF. Considering the first part of the statement, we set $Y = H(\cdot, X)$. By assumption, $X$ is a *strong cubic variation* process. Since $\sigma$ is of class $C^{1,2}$, $H$ is in $C^{1,3}$, and so by applying the Itô formula for *strong cubic variation* processes (see Proposition 2.9), property 4 of Definition 4.1 and the decomposition of the symmetric integral into a classical integral and a covariation term (see Remark 2.7.2), we deduce the following expression for $Y$:

$$Y = H(0, \eta) + \xi + \int_0^{\cdot} \beta(s, X_s) \, dM_s + \tfrac{1}{2}[\beta(\cdot, X), M] + \int_0^{\cdot} \alpha(s, X_s) \, dV_s$$

$$+ \int_0^{\cdot} \partial_s H(s, X_s) \, ds - \frac{1}{4} \int_0^{\cdot} (\sigma^2 (\partial_x \sigma)(\partial_x^{(2)} H))(s, X_s) \, d[\xi, \xi, \xi]_s$$

$$- \frac{1}{12} \int_0^{\cdot} (\partial_x^{(3)} H(s, X_s)) \, d[X, X, X]_s.$$

By property 3, $Y$ is a *strong cubic variation* process as sum of a *strong cubic variation* process and of an $\mathbb{F}$-semimartingale. Moreover, by Remarks 2.5.1 and 2.5.4, $[Y, Y, Y] = [\xi, \xi, \xi]$. Proposition 3.4 tells us that

$$[X, X, X] = [K(\cdot, Y), K(\cdot, Y), K(\cdot, Y)] = \int_0^{\cdot} (\partial_y K(s, Y_s))^3 \, d[Y, Y, Y]_s$$

$$= \int_0^{\cdot} (\sigma(s, X_s))^3 \, d[\xi, \xi, \xi]_s.$$

Using previous equality and computing the partial derivative of $H$ with respect to $x$, we finally reach the result.  □

Before dealing with the case of a possibly vanishing diffusion coefficient $\sigma$, we state the lemma below which will be useful for it.

LEMMA 4.5.  *Let $(X_t, 0 \leq t \leq 1)$ be a solution of equation* (21) *on the probability space* $(\Omega, \mathcal{F}, P)$. *Let $B \in \mathcal{F}$, and $\tau$ a random time. Then, according to the notation of Section 2, the following statements are true:*

1. *the process $X^B$ fulfills properties 2–4 of Definition 4.1 with respect to $\xi^B$, $M^B$ and $V^B$ on the space $(B, \mathcal{F}^B, P^B)$;*
2. *the process $X^\tau$ fulfills properties 2–4 of Definition 4.1 with respect to $\xi^\tau$, $M^\tau$, and $V^\tau$;*
3. *if the coefficients of equation* (21) *are autonomous, and $X$ fulfills property 4 only for autonomous functions, then the process $X_{\cdot + \tau}$ fulfills properties 2 and 3 of Definition 4.1, and property 4 only for autonomous functions, with respect to the processes $\xi_{\cdot + \tau}$, $M_{\cdot + \tau} - M_\tau$ and $V_{\cdot + \tau}$.*

PROOF.  The first and the last point are direct consequences of Lemmas 2.11, 2.12 and 2.13. Concerning the second one, we clearly have that $X^\tau$ is a *strong cubic variation* process by Lemma 2.13. By Lemma 2.12,

$$[\beta(\cdot, X), M]^\tau = [\beta^\tau, M^\tau],$$

with $\beta^\tau = \beta(\cdot \wedge \tau, X_{\cdot \wedge \tau})$. Moreover, the continuity of $M$ and $\beta$ ensures the convergence to zero, almost surely, of the sequence of processes

$$\frac{1}{\varepsilon} \int_0^{\cdot} (\beta((s + \varepsilon) \wedge \tau, X_{(s+\varepsilon) \wedge \tau}) - \beta(s \wedge \tau, X_{s \wedge \tau}))(M_{(s+\varepsilon) \wedge \tau} - M_{s \wedge \tau}) \, ds$$

$$- \frac{1}{\varepsilon} \int_0^{\cdot} (\beta(s + \varepsilon, X_{(s+\varepsilon) \wedge \tau}) - \beta(s, X_{s \wedge \tau}))(M_{(s+\varepsilon) \wedge \tau} - M_{s \wedge \tau}) \, ds$$

$$= \frac{1}{\varepsilon} \int_{(\tau - \varepsilon) \wedge \cdot}^{\tau \wedge \cdot} (\beta(\tau, X_\tau) - \beta(s + \varepsilon, X_\tau))(M_\tau - M_s) \, ds.$$

This implies that $[\beta(\cdot, X^\tau), M^\tau] = [\beta^\tau, M^\tau] = [\beta(\cdot, X), M]^\tau$ exists and it has bounded variation.

If $\psi$ is in $C^{1,\infty}([0,1] \times \mathbb{R})$, at the same way we have

$$\frac{1}{\varepsilon} \int_0^\cdot (\psi(s \wedge \tau, X_{s \wedge \tau}) - \psi(s, X_{s \wedge \tau}))(X_{(s+\varepsilon) \wedge \tau} - X_{(s-\varepsilon) \wedge \tau}) \, ds$$

$$= \frac{1}{\varepsilon} \int_{\tau \wedge \cdot}^{(\tau+\varepsilon) \wedge \cdot} (\psi(\tau, X_\tau) - \psi(s, X_\tau))(X_\tau - X_{s-\varepsilon}) \, ds,$$

and the right-hand side of the equality converges uniformly to zero almost surely, then

$$\int_0^\cdot \psi(s, X_s^\tau) \, d^\circ X_s^\tau = \left( \int_0^\cdot \psi(s, X_s) \, d^\circ X_s \right)^\tau,$$

and so using successively Lemmas 2.11 and 2.12, we obtain that $X^\tau$ fulfills also property 4 of Definition 4.1. $\square$

To treat the case when $\sigma$ is possibly vanishing, we define

$$\nu^\sigma := I_{\{\eta \in S\}}(\omega) H(0, \eta) \qquad \text{for every } \omega \text{ in } \Omega.$$

PROPOSITION 4.6. *Let $\sigma$ be in $C^{1,2}$ satisfying assumptions (H$_1$) and (H$_2$), and $\beta$ be in $C^{0,1}$. Then if $(X_t, 0 \leq t \leq 1)$ is a solution to equation* (21), *adapted to $\mathbb{F}$, and $P(\{\eta \in S^n\}) = 1$, for some $n \geq 0$, it holds*

$$P(\{X_t \in S^n, \forall t \in [0,1]\}) = 1$$

*and*

$$H(\cdot, X) = \xi + N \qquad \text{for all } t \text{ in } [0,1], \ a.s.,$$

*where $N$ is the $\mathbb{F}$-semimartingale*

$$N = \nu^\sigma + \int_0^\cdot \beta(s, X_s) \, dM_s + \int_0^\cdot \alpha(s, X_s) \, dV_s + \int_0^\cdot \partial_s H(s, X_s) \, ds$$

$$+ \tfrac{1}{2}[\beta(\cdot, X), M] + \tfrac{1}{12} \int_0^\cdot (\sigma \partial_x^{(2)} \sigma + (\partial_x \sigma)^2)(s, X_s) \, d[\xi, \xi, \xi]_s.$$

*Furthermore, if $\sigma$ is autonomous, the result still holds even if $X$ fulfills property 4 of Definition 4.1 only for autonomous functions.*

PROOF. Let $D = \mathbb{R}/S$. For every $h$ in $\mathbb{N}^*$, let $\tau^h$ be the first instant the distance between the process $X$ and $D$ becomes smaller than $h^{-1}$:

$$\tau^h = \inf\{t \in [0,1], \text{ s.t. } d(X_t, D) \leq h^{-1}\} \wedge 1,$$

where for every $C$ closed set of $\mathbb{R}$, $x \mapsto d(x, C) = \inf_{r \in C} |x - r|$ is continuous and its support is equal to $C$. We denote, according to the notations of Section 2, $\Omega^h = \{\tau^h > 0\}$, $\mathcal{F}_t^h = \mathcal{F}_t^{\Omega^h}$, $\mathbb{F}^h = (\mathcal{F}_t^h)_{0 \le t \le 1}$, $P^h = P^{\Omega^h}$, and for every stochastic process $Y$ on $\Omega$, we put $Y^h = (Y^{\Omega^h})^{\tau^h}$. Since $P(\eta \in S) = 1$, there exists $k > 0$ such that $P(\Omega^h) > 0$, for every $h \ge k$.

Let $h \ge k \vee (d(c_n, D))^{-1}$ be fixed. We observe that $\Omega^h$ is $\mathcal{F}_0$-measurable; hence, $\mathbb{F}^h$ belongs to $\mathcal{S}(M^h)$. Suppose that $X$ is a solution to equation (21). By Lemmas 4.5.1 and 4.5.2, $X^h$ is a solution of

$$
\begin{cases}
d^\circ X_t^h = \sigma(t, X_t^h)[d^\circ \xi_t^h + \beta(t, X_t^h)\, d^\circ M_t^h + \alpha(t, X_t^h)\, dV_t^h], & 0 \le t \le 1, \\
X_0^h = \eta^h,
\end{cases}
$$

on the probability space $(\Omega^h, \mathcal{F}^h, P^h)$. Moreover, by construction,

$$
P^h(\{X_t^h \in S^{n,h}, \forall t \in [0,1]\}) = 1,
$$

with $S^{n,h} = \{x \in S^n, \text{ s.t. } d(x, D) \ge h^{-1}\}$. Let $\sigma^h : [0,1] \times \mathbb{R} \to \mathbb{R}$ be a function with the same regularity as $\sigma$, never vanishing, and agreeing with $\sigma$ on $S^{n,h}$ together its first and second derivatives in $x$, and its first derivative in $t$. Then $X^h$ is still a solution of

$$
\begin{cases}
d^\circ X_t^h = \sigma^h(t, X_t^h)[d^\circ \xi_t^h + \beta(t, X_t^h)\, d^\circ M_t^h + \alpha(t, X_t^h)\, dV_t^h], & 0 \le t \le 1, \\
X_0^h = \eta^h.
\end{cases}
$$

If $X$ fulfills property 4 only for autonomous functions, then, by Lemma 2.12, $X^h$ carries on doing it with respect to the processes $\xi^h$, $M^h$ and $V^h$, even after having replaced $\sigma$ by $\sigma^h$. In particular, Lemma 4.4 can be applied in both of these two cases. Consequently, if

$$
H^{n,h}(t, x) = \int_{c_n}^x \frac{1}{\sigma^h(t, z)}\, dz,
$$

on $\Omega^h$ it holds $P^h$ almost surely:

$$
\begin{aligned}
H^{n,h}(\cdot, X^h) = {} & H^{n,h}(0, \eta^h) + \xi^h + \int_0^\cdot \beta(s, X_s^h)\, dM_s^h + \int_0^\cdot \alpha(s, X_s^h)\, dV_s^h \\
& + \int_0^\cdot \partial_s H^{n,h}(s, X_s^h)\, ds + \tfrac{1}{2}[\beta(\cdot, X^h), M^h] \\
& + \tfrac{1}{12} \int_0^\cdot (\sigma^h \partial_x^{(2)} \sigma^h + (\partial_x \sigma^h)^2)(s, X_s^h)\, d[\xi^h, \xi^h, \xi^h]_s.
\end{aligned}
$$

We remark that $\{\tau^h > 0\} \subseteq \{\eta^h \in S^{n,h}\}$, and that $h \ge (d(c_n, D))^{-1}$ implies that $c_n$ belongs to $S^{n,h}$. Furthermore, if $x$ belongs to $S^{n,h}$, then $[c_n, x] \subseteq S^{n,h}$. Therefore, $H^{n,h}(t, x) = H(t, x)$, and $\partial_t H^{n,h}(t, x) = \partial_t H(t, x)$, for every $x$ in $S^{n,h}$. Then using Lemma 2.11, Lemma 2.12, and by similar reasonings to those already used in the proof of Lemma 4.5, we obtain the following equality holding $P^h$ almost surely on $\Omega^h$:

$$
(23) \qquad\qquad H(t, X_t) = \xi_t + N_t, \qquad t \le \tau^h,
$$

with

$$N = \nu^\sigma + \int_0^\cdot \beta(s, X_s)\, dM_s + \int_0^\cdot \alpha(s, X_s)\, dV_s + \int_0^\cdot \partial_s H(s, X_s)\, ds$$

$$+ \tfrac{1}{2}[\beta(\cdot, X), M] + \tfrac{1}{12} \int_0^\cdot (\sigma \partial_x^{(2)} \sigma + (\partial_x \sigma)^2)(s, X_s)\, d[\xi, \xi, \xi]_s.$$

Let $\tau = \lim_{h \to +\infty} \tau^h$. Since $\bigcup_{h=0}^{+\infty} \Omega^h = \Omega$, almost surely, we get, for $t = \tau^h$,

$$\lim_{h \to +\infty} H(\tau^h, X_{\tau^h}) = \xi_\tau + N_\tau \qquad \text{a.s.}$$

On the other hand, thanks to the continuity of $X$, $d(X_\tau, D) = 0$ on $\{\tau < 1\}$. This implies

$$\{\tau < 1\} \cup (\{\tau = 1\} \cap \{X_\tau \in D\}) \subseteq \{X_\tau \in \partial D\}.$$

Furthermore, by assumption ($H_2$),

$$\{X_\tau \in \partial D\} \subseteq \left\{ \lim_{h \to +\infty} |H(\tau^h, X_{\tau^h})| = +\infty \right\}$$

$$\subseteq \left\{ \lim_{h \to +\infty} H(\tau^h, X_{\tau^h}) = \xi_\tau + N_\tau \right\}^c.$$

Then it must hold $P(\{\tau < 1\} \cup (\{\tau = 1\} \cap \{X_1 \in D\})) = 0$. We thus have obtained the first part of our result since

$$(\{\tau < 1\} \cup (\{\tau = 1\} \cap \{X_1 \in D\}))^c = \{X_t \in S^n, \forall t \in [0, 1]\}.$$

To complete the proof, it is sufficient to take the limit for $h \to +\infty$ in (23). $\square$

PROPOSITION 4.7. *Let $\sigma$, $\alpha$ and $\beta$ be autonomous, $\sigma$ in $C^{1,2}$, satisfying assumption ($H_2$), and $\beta$ in $C^{0,1}$. Let $X$ be a solution to (21) adapted to $\mathbb{F}$. Then if $P(\{\eta \in D\}) = 1$,*

$$P(\{X_t \in D, \ \forall t \in [0, 1]\}) = 1,$$

*and so $X_t = \eta, \forall t \in [0, 1]$, almost surely.*

PROOF. For every $h \in \mathbb{N}^*$, we consider the first instant the distance between the process $X$ and $D$ becomes greater than $h^{-1}$:

$$\tau^h = \inf\{t \in [0, 1] \text{ s.t. } d(X_t, D) \geq h^{-1}\} \wedge 1,$$

and we put $Y_t^h = Y_{t + \tau^h}$, for $Y = X, \xi, V$ and $M_t^h = M_{t + \tau^h} - M_{\tau^h}$. We observe that $X^h$ is adapted to $\mathbb{F}^h = (\mathcal{F}_t^h)_{t \in [0, 1]}$, where

$$\mathcal{F}_t^h = \{A \in \mathcal{F} | A \cap \{\tau^h \leq s - t\} \in \mathcal{F}_s, \forall s \geq t\},$$

and that $\mathbb{F}^h$ belongs to $\mathcal{S}(M^h)$ (see Problem 3.27 of [19]). Then combining Lemma 4.5.3 and Proposition 4.6, we find that

$$P(\{X_{\tau^h} \in S^m\} \cap \{X_t \in S^m, \forall t \geq \tau^h\}) = P(\{X_{\tau^h} \in S^m\}) \qquad \forall h, m \in \mathbb{N}^*.$$

In particular, since $\tau^h \leq \tau^k$ when $h \geq k$,

$$P(\{X_{\tau^h} \in S^m\} \cap \{X_{\tau^k} \in S^n\}) = 0 \qquad \forall n \neq m, h \geq k.$$

This implies

$$(24) \qquad P(\{X_{\tau^k} \in S^n\}) = P\left(\bigcap_{h \geq k}\{X_{\tau^h} \in S^n\}\right) \qquad \forall n \in \mathbb{N}, \ \forall k \in \mathbb{N}^*.$$

Furthermore, again by Proposition 4.6, we get

$$(25) \qquad H(X_1) - H(X_{\tau^h}) - Y^h = 0, \qquad \text{a.s. on } \{X_{\tau^h} \in S^n\}, \ \forall h \in \mathbb{N}^*,$$

with

$$Y^h = \xi_1 + \int_{\tau^h}^1 \beta(X_s)\, dM_s + \int_{\tau^h}^1 \alpha(X_s)\, dV_s + \tfrac{1}{2}([\beta(X), M]_1 - [\beta(X), M]_{\tau^h})$$
$$+ \tfrac{1}{12}\int_{\tau^h}^1 (\sigma\partial_x^{(2)}\sigma + (\partial_x\sigma)^2)(X_s)\, d[\xi, \xi, \xi]_s.$$

Using assumption (H$_2$) and equality (24), we thus find

$$P(\{X_{\tau^k} \in S^n\}) = P\left(\bigcap_{h \geq k}\{X_{\tau^h} \in S^n\}\right) = 0 \qquad \forall k, n \in \mathbb{N},$$

since in the subspace $\bigcap_{h \geq k}\{X_{\tau^h} \in S^n\}$ we are allowed to take the limit in equality (25). This holds for every $k$ and $n$ in $\mathbb{N}^*$, so we get

$$P(\{X_t \in D, \forall t \in [0,1]\}^c) \leq P\left(\bigcup_{k > 0}\{X_{\tau^k} \in S\}\right) = 0. \qquad \square$$

### 4.4. *Existence and uniqueness.*

PROPOSITION 4.8. *Suppose that there exists a filtration $\mathbb{H} \supseteq \mathbb{F}$, with respect to which the vector $(\xi, M)$ satisfies the hypothesis* (D). *Let $\sigma$ be in $C^{1,2}$, satisfying assumptions* (H$_1$) *and* (H$_2$)*, $\beta$ be in $C^{0,1}$. If $(Y_t, 0 \leq t \leq 1)$ is an $\mathbb{F}$-adapted solution of the stochastic differential equation*

$$Y = \nu^\sigma + \xi + \int_0^{\cdot} \widetilde{\beta}(s, Y_s)\, dM_s + \int_0^{\cdot} \widetilde{\alpha}(s, Y_s)\, dV_s + \int_0^{\cdot} \widetilde{\partial_s H}(s, Y_s)\, ds$$

$$(26) \qquad + \tfrac{1}{2}\int_0^{\cdot} \widetilde{\partial_x\beta}\widetilde{\beta}\widetilde{\sigma}(s, Y_s)\, d[M, M]_s + \tfrac{1}{2}\int_0^{\cdot} \widetilde{\partial_x\beta}\widetilde{\sigma}(s, Y_s)\, d[M, \xi]_s$$

$$+ \tfrac{1}{12}\int_0^{\cdot} (\widetilde{\sigma}\widetilde{\partial_x^{(2)}\sigma} + \widetilde{(\partial_x\sigma)}^2)(s, Y_s)\, d[\xi, \xi, \xi]_s,$$

*then the process*

$$(27) \qquad X = \sum_{n=0}^{\infty} I_{\{\eta \in S^n\}} K^n(\cdot, Y) + I_{\{\eta \in D\}} \eta$$

*is a solution of equation* (21) *adapted to* $\mathbb{F}$*; Conversely, if* $P(\{\eta \in S\}) = 1$ *or* $\sigma$, $\beta$ *and* $\alpha$ *are autonomous and* $(X_t, 0 \le t \le 1)$ *is a solution to equation* (21)*, adapted to* $\mathbb{F}$*, then the process*

$$Y = I_{\{\eta \in S\}} H(\cdot, X) + I_{\{\eta \in D\}} \xi$$

*solves equation* (26)*, and it is* $\mathbb{F}$*-adapted.*

PROOF. Let $(Y_t, 0 \le t \le 1)$ be an $\mathbb{F}$-adapted solution of equation (26). Define $(X_t, 0 \le t \le 1)$ as in formula (27). $X$ is a continuous process with $X_0 = \eta$. Furthermore, $Y$ is a *strong cubic variation* process as the sum of $\xi$ and a semimartingale (recall Remark 2.5.1) and so, by Proposition 3.4, the process $K^n(\cdot, Y)$, for every $n$, has a finite strong cubic variation too. Then $X$ has the same property and

$$[X, X, X] = \sum_{n=0}^{\infty} I_{\{\eta \in S^n\}} [K^n(\cdot, Y), K^n(\cdot, Y), K^n(\cdot, Y)]$$

$$= \sum_{n=0}^{\infty} I_{\{\eta \in S^n\}} \int_0^{\cdot} (\sigma(s, K^n(s, Y_s))^3 \, d[\xi, \xi, \xi]_s$$

$$= \int_0^{\cdot} (\sigma(s, X_s))^3 [\xi, \xi, \xi]_s,$$

where for the last equality we used the fact that $\sigma(t, X_t) I_{\{\eta \in D\}} = 0$, for every $0 \le t \le 1$. Thanks to hypothesis ($\mathcal{D}$), $Y$ is the sum of $R$ and the process $\widetilde{Q} = Y - R$, with $\widetilde{Q} = Q + \int_0^{\cdot} h_s \, dM_s + \widetilde{V}$, $h$ continuous and $\mathbb{H}$-adapted, and $\widetilde{V}$ having bounded variation. Proposition 3.9 implies that $(\widetilde{Q}, M)$ has all its mutual brackets. Then the vector $(Y, M)$ verifies the hypothesis ($\mathcal{D}$), with respect to $\mathbb{H}$. By Proposition 3.18, $[\beta(\cdot, X), M]$ has bounded variation since it is equal to

$$\sum_{n=0}^{\infty} I_{\{\eta \in S^n\}} [\beta(\cdot, K^n(\cdot, Y)), M],$$

with

$$[\beta(\cdot, K^n(\cdot, Y)), M] = \int_0^{\cdot} (\partial_x \beta \sigma)(s, K^n(s, Y_s)) \, d[Y, M]_s$$

$$(28) \qquad = \int_0^{\cdot} (\partial_x \beta \sigma)(s, K^n(s, Y_s)) \, d[\xi, M]_s$$

$$+ \int_0^{\cdot} (\beta \partial_x \beta \sigma)(s, K^n(s, Y_s)) \, d[M, M]_s,$$

on $\{\eta \in S^n\}$. Let $\psi$ be of class $C^{1,\infty}$. We first remark that, since both classical and symmetric integral have a local character (see [23] for the classical integral and Lemma 2.11 for the symmetric one), for every $n$ in $\mathbb{N}^*$ on $\{\eta \in S^n\}$, it holds:

$$Y = \nu^\sigma + \xi + \int_0^{\cdot} \beta(s, K^n(s, Y_s)) \, dM_s + \int_0^{\cdot} \alpha(s, K^n(s, Y_s)) \, dV_s$$

$$+ \int_0^{\cdot} \partial_s H(s, K^n(s, Y_s)) \, ds + \frac{1}{2} \int_0^{\cdot} \partial_x \beta \beta \sigma(s, K^n(s, Y_s)) \, d[M, M]_s$$

$$+ \frac{1}{2} \int_0^{\cdot} \partial_x \beta \sigma(s, K^n(s, Y_s)) \, d[M, \xi]_s$$

$$+ \frac{1}{12} \int_0^{\cdot} (\sigma \partial_x^{(2)} \sigma + (\partial_x \sigma)^2)(s, K^n(s, Y_s)) \, d[\xi, \xi, \xi]_s.$$

We apply the Itô formula for *strong cubic variation* processes to write

$$X = \sum_{n=0}^{\infty} I_{\{\eta \in S^n\}} X^n + I_{\{\eta \in D\}} \eta,$$

with

$$X^n = \eta + \int_0^{\cdot} \partial_s K^n(s, Y_s) \, ds + \int_0^{\cdot} \partial_y K^n(s, Y_s) \, d^\circ Y_s$$

$$- \frac{1}{12} \int_0^{\cdot} \partial_y^{(3)} K^n(s, Y_s) \, d[Y, Y, Y]_s.$$

Using equality (28), we can write on $\{\eta \in S^n\}$,

$$Y = \nu^\sigma + \xi + \int_0^{\cdot} \beta(s, K^n(s, Y_s)) \, d^\circ M_s + \int_0^{\cdot} \alpha(s, K^n(s, Y_s)) \, dV_s$$

$$+ \int_0^{\cdot} \partial_s H(s, K^n(s, Y_s)) \, ds$$

$$+ \frac{1}{12} \int_0^{\cdot} (\sigma \partial_x^{(2)} \sigma + (\partial_x \sigma)^2)(s, K^n(s, Y_s)) \, d[\xi, \xi, \xi]_s.$$

Deriving with respect to $s$ the equality $H(s, K^n(s, y)) = y$, we obtain the relation

$$\partial_s K^n(s, y) = -\sigma(s, K^n(s, y)) \partial_s H(s, K^n(s, y)),$$

which combined with equation (26), the equalities

$$\partial_y K^n(s, y) = \sigma(s, K^n(s, y)), \qquad \partial_y^{(2)}(\sigma(s, K^n(s, y))) = \partial_y^{(3)} K^n(s, y),$$

and Corollary 3.19 gives the following expression for $X^n$ on $\{\eta \in S^n\}$:

$$X^n = \eta + \int_0^{\cdot} \sigma(s, X_s^n) \, d^\circ \xi_s + \int_0^{\cdot} (\sigma \beta)(s, X_s^n) \, d^\circ M_s + \int_0^{\cdot} (\sigma \alpha)(s, X_s^n) \, dV_s.$$

Coefficients appearing in the last expression for $X^n$ and function $\psi$ are regular enough to use successively Lemma 2.10 and Corollary 3.19 to get on $\{\eta \in S^n\}$:

$$\int_0^\cdot \psi(t, X_t^n) \, d^\circ X_t^n = \int_0^\cdot (\psi\sigma)(t, X_t^n)[d^\circ \xi_t + \beta(t, X_t^n) \, d^\circ M_t + \alpha(t, X_t^n) \, dV_t]$$
$$- \tfrac{1}{4} \int_0^\cdot (\partial_x \sigma)(\sigma^2)(\partial_x \psi)(t, X_t^n) \, d[\xi, \xi, \xi]_t.$$

The conclusion follows since $\int_0^\cdot \psi(t, X_t) \, d^\circ X_t = \sum_{n=0}^\infty I_{\{\eta \in S^n\}} \int_0^\cdot \psi(t, X_t^n) \, d^\circ X_t^n$, almost surely.

We consider the second part of the statement. By Proposition 4.6,

$$Y = H(\cdot, X) = \xi + N \qquad \text{on } \{\eta \in S\}.$$

The vector $(\xi, N, M)$ fulfills the hypothesis $(\mathcal{D})$ with respect to $\mathbb{H}$. Indeed, $N = \int_0^\cdot h_s \, dM_s + \widetilde{V}$, with $h$ continuous and $\mathbb{H}$-adapted, and $\widetilde{V}$ with bounded variation. By Proposition 3.18,

$$I_{\{\eta \in S\}}[\beta(\cdot, X), M] = \sum_{n=0}^\infty I_{\{\eta \in S^n\}}[\beta(\cdot, K^n(\cdot, \xi + N)), M],$$

with

$$[\beta(\cdot, K^n(\cdot, \xi + N)), M] = \int_0^\cdot (\partial_x \beta \sigma(s, K^n(s, \xi_s + N_s))) \, d[\xi, M]_s$$
$$+ \int_0^\cdot (\beta \partial_x \beta \sigma(s, K^n(s, \xi_s + N_s))) \, d[M, M]_s.$$

Therefore, on $\{\eta \in S^n\}$, $N$ is more explicitly given by the following expression:

$$N = \nu^\sigma + \int_0^\cdot \beta(s, K^n(s, Y_s)) \, dM_s + \int_0^\cdot \alpha(s, K^n(s, Y_s)) \, dV_s$$
$$+ \tfrac{1}{2} \int_0^\cdot (\sigma\beta\partial_x\beta)(s, K^n(s, Y_s)) \, d[M, M]_s$$
$$\text{(29)} \qquad + \tfrac{1}{2} \int_0^\cdot (\sigma\partial_x\beta)(s, K^n(s, Y_s)) \, d[M, \xi]_s$$
$$+ \int_0^\cdot \partial_s H(s, K^n(s, Y_s)) \, ds$$
$$+ \tfrac{1}{12} \int_0^\cdot (\sigma\partial_x^{(2)}\sigma + (\partial_x\sigma)^2)(s, K^n(s, Y_s)) \, d[\xi, \xi, \xi]_s.$$

Putting expression (29) in the equality

$$Y = I_{\{\eta \in S\}}(\xi + N) + I_{\{\eta \in D\}}\xi = \xi + \sum_{n=0}^\infty I_{\{\eta \in S^n\}} N,$$

we achieve the proof of the proposition. □

Theorem 4.9. *Suppose that there exists a filtration $\mathbb{H} \supseteq \mathbb{F}$, with respect to which the vector $(\xi, M)$ satisfies the hypothesis (D). Let $\sigma$ satisfy assumptions (H$_1$), (H$_2$), and the following hypotheses:*

(H$_3$)
  (i)   $\sigma$ is in $C^{1,2}$,
  (ii)  $\partial_x^{(2)} \sigma$ is locally Lipschitz in $x$, uniformly in $t$ ,
  (iii) $\sup_{(t,x)\in[0,1]\times S^n} |\partial_t \log(|\sigma(t,x)|)| \leq a^n, \forall n \in \mathbb{N}$,
  (iv)  $(|\partial_x \sigma|^2 + |\sigma \partial_x^{(2)} \sigma|)(t,x) \leq a_n(1 + |H^n(t,x)|)$,
        $(t,x) \in [0,1] \times S^n, \ \forall n \in \mathbb{N}$,

*for some sequences $(a_n)_{n\in\mathbb{N}}$, in $\mathbb{N}$; let $\beta$ and $\alpha$ verify the following:*

(H$_4$)
  (i)   $\beta$ is in $C^{0,1}$ and it is bounded,
  (ii)  $\partial_x \beta$ and $\alpha$ are locally Lipschitz in $x$, uniformly in $t$,
  (iii) $(|\sigma||\partial_x \beta| + |\alpha|)(t,x) \leq a_n(1 + |H^n(t,x)|), (t,x) \in [0,1] \times S^n$,

*for all $n$ in $\mathbb{N}$. Then if $P(\{\eta \in S\}) = 1$ or $\sigma$, $\beta$ and $\alpha$ are autonomous, equation (21) has a unique $\mathbb{F}$-adapted solution given by*

$$X = \sum_{n=0}^{\infty} I_{\{\eta\in S^n\}} K^n(\cdot, Y) + I_{\{\eta\in D\}}\eta,$$

*where $Y$ is the unique $\mathbb{F}$-adapted solution to equation (26).*

Remark 4.10. We emphasize that hypothesis (H$_4$) has to be satisfied by $\alpha$ a.s. In the sequel we will implicitly use this convention.

Proof of Theorem 4.9. The result follows from the existence and uniqueness of equation (26). The last holds since assumptions (H$_3$) and (H$_4$) imply the local Lipschitz continuity and the linear growth property of the coefficients of equation (26), which are sufficient conditions to ensure its existence and uniqueness (see [10], page 29, Lemma 34). In fact, the functions

$$(t,y) \mapsto \beta(t, K^n(t,y)),$$

$$\sigma \partial_x^{(2)} \sigma(t, K^n(t,y)), \alpha(t, K^n(t,y)), (\partial_x \sigma(t, K^n(t,y)))^2$$

and

$$(t,y) \mapsto \sigma \partial_x \beta(t, K^n(t,y))$$

have linear growth thanks to the boundedness of $\beta$, (iv) of (H$_3$) and (iii) of (H$_4$); moreover, they are locally Lipschitz, being the composition of continuous functions differentiable with continuity or locally Lipschitz in $y$. The map $(t, y) \mapsto \partial_t H^n(t, K^n(t, y))$ is locally Lipschitz, being differentiable with continuity with respect to $y$. By (iii) of (H$_3$), $|\partial_t H^n(t, x)| \leq a_n |H^n(t, x)|$, which implies the linear growth for $(t, y) \mapsto \partial_t H(t, K^n(t, y))$. $\square$

Recalling Examples 3.8 and 3.7, one can prove the following results.

COROLLARY 4.11.  *Suppose that there exist two adapted processes $Q$ and $R$, such that $\xi = R + Q$, $R$ is $\mathcal{F}_0$-measurable and $(Q, M)$ has all its mutual brackets. Let $\sigma$, $\beta$ and $\alpha$ verify the regularity assumptions of Proposition 4.9. Then if the $P(\{\eta \in S\}) = 1$ or the coefficients are autonomous, there exists a unique $\mathbb{F}$-adapted solution to equation* (21).

COROLLARY 4.12.  *Suppose that there exist two adapted processes $Q$ and $R$, such that $\xi = R + Q$, with $R$ independent from $M$, $(Q, M)$ having all its mutual brackets, and $\mathbb{F} \subseteq \mathbb{H}$, being $\mathcal{H}_t = \sigma(M_s, 0 \leq s \leq t) \vee \sigma(R)$, for every $0 \leq t \leq 1$. Let $\sigma$, $\beta$ and $\alpha$ verify the regularity assumptions of Proposition 4.9. Then if the $P(\{\eta \in S\}) = 1$ or the coefficients are autonomous, there exists a unique $\mathbb{F}$-adapted solution to equation* (21).

If $\sigma$ is bounded from below from a positive constant, we can solve with our methods an equation already studied in [7], where the diffusion coefficient does not appear as a multiplier factor. There the coefficient $\beta$ was equal to zero, $\sigma$ autonomous and of class $C^{1,3}$. The authors needed to introduce the notion of  *strong cubic vector Itô processes* in the Definition 4.1, requiring more than the finite cubic variation of a solution $X$. In particular, existence and uniqueness were proved to hold in a smaller class than ours, with more regularity on $\sigma$.

4.5. *On the uniqueness of the integral equation.* We aim here at adding hypotheses on the coefficients driving equation (21) to find a suitable class of processes among which its solution, in the sense described in Definition 4.1, exists and it is the unique solution to the integral equation (22).

REMARK 4.13.  1. Let $Z$ be in $\mathcal{C}_\xi^2(\mathbb{H})$ and $\psi$ in $C^{1,4}$, with $\partial_t \psi$ in $C^{0,2}$. Then the process $(\psi(t, Z_t), 0 \leq t \leq 1)$ is in $\mathcal{C}_\xi^2(\mathbb{H})$.
2. Let

$$(X^k(t, x), 0 \leq t \leq 1, x \in \mathbb{R})_{k \in \mathbb{N}}$$

be a sequence of $C^2$ $\mathbb{H}$-Itô-semimartingale fields, of this form,

$$X^k(t, x) = f^k(x) + \sum_{i=1}^n \int_0^t a^{k,i}(s, x) \, dN_s^i + \sum_{j=1}^m \int_0^t b^{k,j}(s, x) \, dV_s^j,$$

and $(\Omega_k)_{k \in \mathbb{N}}$ be a sequence of subspaces of $\Omega$ in $\mathcal{H}_0$, with $\bigcup_{k=0}^{\infty} \Omega_k = \Omega$, a.s. Then the random field

$$Y(t,x) = \sum_{k=0}^{\infty} I_{\Omega_k} X^k(t,x)$$

is a $C^2$ $\mathbb{H}$-Itô-semimartingale field of the form

$$Y(t,x) = f(x) + \sum_{i=1}^{n} \int_0^t a^i(s,x) \, dN_s^i + \sum_{j=1}^{m} \int_0^t b^j(s,x) \, dV_s^j$$

with

$$f(x) = \sum_{k=0}^{\infty} I_{\Omega_k} f^k(x),$$

$$a^i(t,x) = \sum_{k=0}^{\infty} I_{\Omega_k} a^{k,i}(t,x),$$

$$b^j(t,x) = \sum_{k=0}^{\infty} I_{\Omega_k} b^{k,j}(t,x).$$

PROPOSITION 4.14. *Suppose that there exists a filtration $\mathbb{H} \supseteq \mathbb{F}$, with respect to which the vector $(\xi, M)$ satisfies the hypothesis $(\mathcal{D})$. Let $\sigma$, $\beta$ and $\alpha$ satisfy hypotheses of Proposition 4.9, with, furthermore, $\sigma$ in $C^{1,4}$, $\partial_t \sigma$ in $C^{0,2}$, $\beta$ in $C^{1,3}$ and $\partial_t \beta$ in $C^{0,1}$. Then if $P(\{\eta \in S\}) = 1$, or $\alpha$, $\sigma$ and $\beta$ are autonomous, there exists a unique $\mathbb{F}$-adapted solution to the integral equation (22) in the space $\mathcal{C}^2_{\xi, \eta}(\mathbb{H})$ of all processes in $\mathcal{C}^2_\xi(\mathbb{H})$, starting at $\eta$.*

PROOF. The existence was proved in Proposition 4.9. Consider, in fact, the process $Y$ which is the unique solution of equation (26). Classical Itô formula for semimartingales applied to the function $K^n$ and the semimartingale $N = Y - \xi$, shows that the random field $(K^n(t, x + N_t), t \in [0,1], x \in \mathbb{R})$ is a $C^2$ $\mathbb{H}$-Itô-semimartingale field driven by the local martingale $M$. Therefore, by Remark 4.13, $X$ is in $\mathcal{C}^2_{\xi, \eta}(\mathbb{H})$.

Regarding uniqueness, we show that an integral solution in $\mathcal{C}^2_{\xi, \eta}(\mathbb{H})$ is a solution in the sense described in Definition 4.1. Let $Z$ be the random field in $\mathcal{C}^2_{\xi, \eta}(\mathbb{H})$ such that $X = Z(\cdot, \xi)$, where $X$ is a solution to equation (22). Condition 1 is fulfilled by the hypothesis. Since $\xi$ is an $\mathbb{H}$-adapted *strong cubic variation* process and $Z$ is a $C^2$ $\mathbb{H}$-Itô-semimartingale field, by Proposition 3.4, $X$ satisfies Condition 2. By the classical Itô formula, $(\beta(t, Z(t,x)), 0 \leq t \leq 1, x \in \mathbb{R})$ is a $C^1$ $\mathbb{H}$-Itô-semimartingale field driven by a vector of local martingales $(N^1, \ldots, N^n)$ such that the vector $(\xi, N^1, \ldots, N^n)$ satisfies the hypothesis $(\mathcal{D})$ with respect to $\mathbb{H}$. By definition, there exist two $\mathbb{H}$-adapted processes $\bar{R}$ and $\bar{Q}$ such that $\xi = \bar{R} + \bar{Q}$, $(\bar{Q}, N^1, \ldots, N^n)$ has all its mutual brackets, and $\bar{R}_{\varepsilon+}$. is $\mathbb{H}$-adapted. By Corollary 3.11, $[\bar{R}, M] = 0$. This implies the existence of $[M, \bar{Q}]$ which equals $[\xi, M]$. Then $(\xi, N^1, \ldots, N^n, M)$ verifies the hypothesis ($\mathcal{D}$) with respect to $\mathbb{H}$, and by Proposition 3.16, Condition 3 is established. Since $\partial_t \sigma$ belongs to $C^{0,2}$, it follows from the details of the proofs that if Condition 4 is fulfilled for functions $\psi$ in $C^{1,\infty}$ with $\partial_t \psi$ in $C^{0,2}$, previous results about uniqueness remain true. Let then $\psi$ be in $C^{1,\infty}$, with $\partial_t \psi$ in $C^{0,2}$. $X$ is a solution of the integral equation, so we can write

$$\int_0^t \psi(s, X_s) \, d^\circ X_s = \int_0^t \widehat{\psi}(s, \xi_s) \, d^\circ \left( \int_0^s \widehat{\sigma}(r, \xi_r) \, d^\circ \xi_r \right)$$
$$+ \int_0^t \widehat{\psi}(s, \xi_s) \, d^\circ \left( \int_0^s \widehat{\beta\sigma}(r, \xi_r) \, d^\circ M_r \right)$$
$$+ \int_0^t \widehat{\psi}(s, \xi_s) \, d^\circ \left( \int_0^s \widehat{\alpha\sigma}(r, \xi_r) \, dV_r \right),$$

with the notation $\widehat{\psi}(t, x) = \psi(t, Z(t, x))$, for every function $\psi : [0, 1] \times \mathbb{R}$. As already remarked before, the processes $(\widehat{\psi}(t, \xi_t), 0 \leq t \leq 1)$, as well as $(\widehat{\sigma}(t, \xi_t), 0 \leq t \leq 1)$, are in $\mathcal{C}^2_{\xi,\eta}(\mathbb{H})$ so as to let us apply Proposition 3.25, at the same way the random field $(\widehat{\beta\sigma}(t, x), 0 \leq t \leq 1, x \in \mathbb{R})$ has the properties needed in Corollary 3.19. Then we obtain

$$\int_0^\cdot \psi(s, X_s) \, d^\circ X_s = \int_0^\cdot (\psi\sigma)(s, X_s) \, d^\circ \xi_s + \int_0^\cdot (\psi\beta\sigma)(s, X_s) \, d^\circ M_s$$
$$+ \int_0^\cdot (\psi\alpha\sigma)(s, X_s) \, dV_s$$
$$- \tfrac{1}{4} \int_0^\cdot (\partial_x \psi)(\partial_x \sigma)(s, X_s)(\partial_x Z(s, \xi_s))^2 \, d[\xi, \xi, \xi]_s.$$

By Proposition 3.4,

$$\int_0^\cdot \partial_x \psi \partial_x \sigma(s, X_s)(\partial_x Z(s, \xi_s))^2 \, d[\xi, \xi, \xi]_s$$
$$= \int_0^\cdot (\partial_x \psi \partial_x \sigma)(s, X_s) \, d[X, X, \xi]_s.$$

Finally, by multi-linearity of the 3-covariation application, and remarks 2.5.1 and 3.23,

$$[X, X, \xi] = \left[ \int_0^\cdot \widehat{\sigma}(s, \xi_s) \, d^\circ \xi_s, \int_0^\cdot \widehat{\sigma}(s, \xi_s) \, d^\circ \xi_s, \xi \right]$$
$$= \int_0^\cdot (\sigma(s, X_s))^2 \, d[\xi, \xi, \xi]_s,$$

and so Condition 4 is proved to hold. This leads to the conclusion of the proof. $\square$

4.6. *The finite quadratic variation case.* In this section we suppose that the vector $(\xi, M)$ has all its mutual brackets, in particular, that $\xi$ is a *finite quadratic variation* process. We observe that, under this assumption, the vector $(\xi, M)$ satisfies the hypothesis ($\mathcal{D}$) with respect to the filtration $\mathbb{H} = \mathbb{F}$. Moreover, $\mathcal{C}_\xi^k(\mathbb{F})$ reduces to the set of all the $C^k$ $\mathbb{F}$-Itô-semimartingale fields driven by a vector of semimartingales $(N^1, \ldots, N^n)$ such that $(\xi, N^1, \ldots, N^n)$ has all its mutual brackets.

Results obtained in the previous section can be improved regarding the regularity required for the *diffusion* coefficient $\sigma$, by using techniques which are similar to those already developed in [27] and [10] about stochastic calculus with respect to finite quadratic variation processes. More precisely, the Itô formula for *finite quadratic variation* processes holds for $C^2$ functions of the space variable, which allows us to reduce of one the degree of regularity of $\sigma$.

DEFINITION 4.15. A continuous stochastic process $(X_t, 0 \leq t \leq 1)$ will be said solution to equation (21) if $X_0 = \eta$, the vector $(X, M)$ has all its mutual brackets, and for every $\psi$ in $C^{1,\infty}$, it holds,

$$\int_0^{\cdot} \psi(s, X_s)\, d^\circ X_s = \int_0^{\cdot} \psi\sigma(s, X_s)[d^\circ \xi_s + \beta(s, X_s)\, d^\circ M_s + \alpha(s, X_s)\, dV_s].$$

REMARK 4.16. Definition 4.1 and 4.15 are equivalent. It is sufficient to use Proposition 3.18, and recall that $[\xi, \xi, \xi] = 0$.

Similarly to the finite cubic variation case, we state the following results.

PROPOSITION 4.17. *Let $\sigma$ be in $C^{1,1}$, satisfying assumptions ($H_1$) and ($H_2$), $\beta$ be in $C^{0,1}$. If $(Y_t, 0 \leq t \leq 1)$ is an $\mathbb{F}$-adapted solution of the stochastic differential equation*

$$
\begin{aligned}
(30) \quad Y &= \nu^\sigma + \xi + \int_0^{\cdot} \widetilde{\beta}(s, Y_s)\, dM_s + \int_0^{\cdot} \widetilde{\alpha}(s, Y_s)\, dV_s + \int_0^{\cdot} \widetilde{\partial_t H}(s, Y_s)\, ds \\
&\quad + \frac{1}{2}\int_0^{\cdot} \widetilde{\partial_x \widetilde{\beta} \widetilde{\sigma}}(s, Y_s)\, d[M, M]_s + \frac{1}{2}\int_0^{\cdot} \widetilde{\partial_x \widetilde{\sigma}}(s, Y_s)\, d[M, \xi]_s,
\end{aligned}
$$

*then the process $X = \sum_{n=0}^{\infty} I_{\{\eta \in S^n\}} K^n(\cdot, Y) + I_{\{\eta \in D\}} \eta$ is a solution of equation (21) adapted to $\mathbb{F}$. Conversely, if $P(\{\eta \in S\}) = 1$ or $\sigma$, $\beta$ and $\alpha$ are autonomous and $(X_t, 0 \leq t \leq 1)$ is a solution to equation (21), adapted to $\mathbb{F}$, then the process $Y = I_{\{\eta \in S\}} H(\cdot, X) + I_{\{\eta \in D\}} \xi$ solves equation (30), and it is $\mathbb{F}$-adapted.*

PROPOSITION 4.18. *Let $\sigma$ be in $C^1$, satisfying assumptions ($H_1$) and ($H_2$), and such that*

$$\sup_{(t,x) \in [0,1] \times S^n} |\partial_t \log(|\sigma(t, x)|)| \leq a^n \qquad \forall n \in \mathbb{N}$$

*for some sequences $(a_n)_{n \in \mathbb{N}}$ in $\mathbb{R}^+$; let $\beta$ and $\alpha$ verify hypothesis* ($\text{H}_4$). *Then if* $P(\{\eta \in S\}) = 1$ *or* $\sigma$, $\beta$ *and* $\alpha$ *are autonomous, equation* (21) *has a unique* $\mathbb{F}$-*adapted solution.*

We aim at comparing the results obtained with our method with those achieved in [10] and [27]. There $\sigma$ was not a multiplier coefficient. Then the comparison can be made if $\sigma$ is bounded from below from a positive constant. In such a case equations studied by those authors are particular cases of equation (21), where the symmetric integral is replaced by the *forward* one; see [25] for definition. Note that a solution $X$ has not to be supposed a finite quadratic variation process by point 4 of Definition 4.1. Indeed, $X$ is a finite quadratic variation process if and only if $\int_0^{\cdot} X_t \, d^- X_t$ exists.

We remember that, for two continuous stochastic processes $X$ and $Y$, if the symmetric integral, $\int_0^{\cdot} X_s \, d^\circ Y_s$, and the forward integral, $\int_0^{\cdot} X_s \, d^- Y_s$, exist, then $\frac{1}{2}[X, Y]$ exists and

$$\int_0^{\cdot} X_s \, d^\circ Y_s = \int_0^{\cdot} X_s \, d^- Y_s + \frac{1}{2}[X, Y].$$

Using this relation, under assumptions of Proposition 4.17, we can state this equivalence between the solution to equation (21) in the symmetric and the forward sense. This notion of solution in Definition 4.15 has to be adapted replacing the symmetric integral with the forward one.

A process $X$ is a solution of equation

$$(31) \qquad \begin{cases} d^- X_t = \sigma(t, X_t)[d^- \xi_t + \beta(t, X_t) \, d^- M_t + \alpha(t, X_t) \, d^- V_t], \\ X_0 = \eta, \end{cases}$$

if and only if it solves

$$(32) \qquad \begin{cases} d^\circ X_t = \sigma(t, X_t)[d^\circ \xi_t + \beta(t, X_t) \, d^\circ M_t + \alpha(t, X_t) \, dV_t] \\ \qquad - \frac{1}{2}\sigma(t, X_t)[\gamma^1(t, X_t) \, dV_t^1 + \gamma^2(t, X_t) \, dV_t^2 + \gamma^3(t, X_t) \, dV_t^3], \\ X_0 = \eta, \end{cases}$$

with $\gamma^1 = \partial_x \sigma, \gamma^2 = 2\partial_x \sigma\beta + \sigma\partial_x \beta, \gamma^3(t, x) = \partial_x \sigma\beta^2 + \sigma\beta\partial_x \beta$, and $V^1 = [\xi, \xi], V^2 = [\xi, M], V^3 = [M, M]$.

This equivalence and Proposition 4.18 imply the following.

REMARK 4.19. Suppose that, besides the hypotheses of Proposition 4.18, $\partial_x \sigma$ is locally Lipschitz in $x$, uniformly in $t$, and

$$|\partial_x \sigma|(t, x) \le a_n(1 + |H^n(t, x)|), \qquad (t, x) \in [0, 1] \times S^n, \forall n \in \mathbb{N}.$$

Then equation (31) has a unique solution. Existence and uniqueness are ensured by equation (32). Moreover, the solution is given by $X = \sum_{n=0}^{\infty} I_{\{\eta \in S^n\}} K^n(\cdot, Y) +$

$I_{\{\eta \in D\}}\eta$, where $(Y_t, 0 \leq t \leq 1)$ is the unique solution of

$$Y = \nu^\sigma + \xi + \int_0^\cdot \widetilde{\beta}(s, Y_s)\, dM_s + \int_0^\cdot \widetilde{\alpha}(s, Y_s)\, dV_s + \int_0^\cdot \widetilde{\partial_s H}(s, Y_s)\, ds$$

$$(33) \qquad - \frac{1}{2} \int_0^\cdot \widetilde{\partial_x \sigma}(s, Y_s)\, d[\xi, \xi]_s - \int_0^\cdot \widetilde{\partial_x \sigma \widetilde{\beta}}(s, Y_s)\, d[M, \xi]_s$$

$$- \frac{1}{2} \int_0^\cdot \widetilde{\partial_x \sigma \widetilde{\beta}^2}(s, Y_s)\, d[M, M]_s.$$

REMARK 4.20. If we assume $\beta$ only continuous, bounded and locally Lipschitz, equation (33) still has a unique solution. Nevertheless, $X$ could fail to solve equation (32); indeed, the bracket $[\beta(\cdot, X), M]$ may not exist under this weaker condition.

In order to avoid this additional conditions on $\beta$, equation (31) has to be studied directly using stochastic calculus with respect to finite quadratic variation processes and forward integrals instead of symmetric ones. By these methods, it is possible to show the following result.

PROPOSITION 4.21. *Suppose that $\sigma$ is in $C^{1,1}$ and it satisfies assumptions (H$_1$) and (H$_2$), that $\beta$ is continuous and bounded, $\beta, \alpha, \partial_x \sigma$ are locally Lipschitz in $x$ uniformly in $t$, and moreover, that*

$$\sup_{(t,x) \in [0,1] \times \mathbb{R}} |\partial_t \log(|\sigma(t,x)|)| < +\infty;$$

$$(|\partial_x \sigma| + |\alpha|)(t,x) \leq a_n(1 + |H^n(t,x)|),$$

$$(t,x) \in [0,1] \times S^n, \ \forall\, n \in \mathbb{N}.$$

*Then equation (31) has a unique solution.*

Moreover, as in the finite cubic variation case, we can also state the following.

PROPOSITION 4.22. *Let $\sigma$, $\beta$ and $\alpha$ satisfy hypotheses of Proposition 4.21, with, furthermore, $\sigma$ in $C^{1,3}$, and $\partial_t \sigma$ in $C^{0,1}$. Then, if $P(\{\eta \in S\}) = 1$, or $\alpha$, $\sigma$ and $\beta$ are autonomous, there exists a unique $\mathbb{F}$-adapted solution to the integral equation*

$$X = \eta + \int_0^\cdot \sigma(s, X_s)\, d^- \xi_s + \int_0^\cdot \sigma\beta(s, X_s)\, d^- M_s + \int_0^\cdot \sigma\alpha(s, X_s)\, dV_s,$$

*in the space $\mathcal{C}^1_{\xi,\eta}(\mathbb{F})$ of all processes in $\mathcal{C}^1_\xi(\mathbb{F})$, starting at $\eta$.*

In [10] the authors show the existence and uniqueness of the integral equation (31), supposing $\sigma$ *autonomous* and in $C^{1,4}$, in the class $\mathcal{C}^2_{\xi,\eta} \subset \mathcal{C}^1_{\xi,\eta}$. In [27] an equation of type (31) is studied with semimartingale coefficient $\beta$ equal to zero and an *autonomous* diffusion coefficient. There $\sigma$ is of class $C^3$, bounded with its partial derivative $\partial_x\sigma$. Moreover, the sense of solution is more restrictive in that it involves the notion of *vector Itô processes* which are not necessary to introduce for the application of our method.

4.7. *The Hölder continuous case.* We intend to apply the methods developed in previous sections to the study of the stochastic differential equation (21) when the processes $\xi$ and $V$ have $\gamma$-Hölder continuous paths, with $\frac{1}{2} < \gamma < 1$, the semimartingale coefficient is equal to zero, and $V_t = t$:

$$(34) \qquad \begin{cases} d^\circ X_t = \sigma(t, X_t)[d^\circ \xi_t + \alpha(t, X_t)\, dt], \\ X_0 = \eta. \end{cases}$$

REMARK 4.23. This method could be extended to the case $V = \int_0^\cdot \psi_s\, ds$, with $\psi \in L^2([0,1])$. Indeed, this would imply $V$ $\gamma$-Hölder continuous with $\gamma > \frac{1}{2}$.

We will see that in this case the use of an Itô formula available for processes having Hölder continuous paths will let us reduce the regularity of $\sigma$. If $0 < \gamma < 1$, $C^\gamma$ will denote the Banach space of all $\gamma$-Hölder continuous functions with the norm

$$\|f\|_\gamma = \sup_{s,t\in[0,1] s\neq t} \frac{|f(t) - f(s)|}{|t - s|^\gamma} + \|f\|_\infty.$$

In this context we will look for existence and uniqueness of integral solutions with $\gamma$-Hölder continuous paths. We first recall some results about integral calculus with respect to Hölder functions contained in [9] and [31].

LEMMA 4.24. *Let $f$ and $g$ be in $C^1$, with $f(0) = 0$, and $\alpha + \gamma > 1$. Then the following inequality holds:*

$$\left| \int_0^t f(r)\, dg(r) \right| \leq C\|f\|_\alpha \|g\|_\gamma t^{1+\varepsilon}$$

*for some positive constant $C$ and $0 < \varepsilon < \alpha + \gamma - 1$.*

COROLLARY 4.25. *Let $f$ and $g$ be in $C^1$, and $\alpha + \gamma > 1$. Then the following inequality holds, for every $t, s$ in $[0,1]$:*

$$\left| \int_s^t f(r)\, dg(r) - f(s)(g(t) - g(s)) \right| \leq C\|f\|_\alpha \|g\|_\gamma |t - s|^{1+\varepsilon}$$

*for some positive constant $C$ and $0 < \varepsilon < \alpha + \gamma - 1$. In particular, $\int_0^\cdot f\, dg$ is a $\gamma$-Hölder function.*

Corollary 4.25 implies the following.

PROPOSITION 4.26. *If $\alpha + \gamma > 1$, the map $F: (f, g) \mapsto \int_0^{\cdot} f \, dg$ defined on $C^1 \times C^1$, with values in $C^{\gamma}$, admits a unique continuous extension to $C^{\alpha} \times C^{\gamma}$.*

PROOF. Let $(f, g)$ and $(h, k)$ in $C^{\alpha} \times C^{\gamma}$. The map $F$ is bilinear, therefore,

$$\|F(f, g) - F(h, k)\|_{\gamma} \leq \|F(f - h, k)\|_{\alpha} + \|F(h, g - k)\|_{\gamma}.$$

Let $s, t$ be in $[0, 1]$. By Corollary 4.25,

$$|F(f - h, g)(t) - F(f - h, g)(s)| \leq C\|f - h\|_{\alpha}\|g\|_{\gamma}|t - s|^{\gamma},$$

and similarly,

$$|F(h, g - k)(t) - F(h, g - k)(s)| \leq C\|h\|_{\alpha}\|g - k\|_{\gamma}|t - s|^{\gamma}.$$

This immediately implies

$$\|F(f - h, g)\|_{\gamma} + \|F(h, g - k)\|_{\gamma}$$
$$\leq 2C(\|g\|_{\gamma} \vee \|h\|_{\alpha})\|(f, g) - (h, k)\|_{C^{\alpha} \times C^{\gamma}}. \qquad \square$$

The unique continuous extension of $F$ will be called the *Young* integral and denoted with $\int_0^{\cdot} f \, d^y g$, for every $f$ in $C^{\alpha}$ and $g$ in $C^{\gamma}$.

REMARK 4.27. If $f$ and $h$ are in $C^{\alpha}$ and $g$ in $C^{\gamma}$, with $\alpha + \gamma > 1$, we have

$$\int_0^{\cdot} f \, d^y \left( \int_0^{\cdot} h \, d^y g \right) = \int_0^{\cdot} fh \, d^y g.$$

The equality holds for $(f, g)$ in $C^1 \times C^1$, and it can be extended to $C^{\alpha} \times C^{\gamma}$ by density arguments.

L. C. Young [30] introduced that integral in a more general setting, that is, for $f, g$ having respectively $p$ and $q$ variation with $p^{-1} + q^{-1} = 1$. It can be proved that the Young integral $\int_0^{\cdot} f \, d^y g$ agrees with the symmetric integral $\int_0^{\cdot} f \, d^{\circ} g$ (see [28]) and that it is a *Riemann–Stieltjes* type integral as specified in the following proposition.

PROPOSITION 4.28. *Let $f$ be in $C^{\alpha}$ and $g$ in $C^{\gamma}$, with $\alpha + \gamma > 1$. Then for every $0 \leq t \leq 1$,*

$$\lim_{\delta \to 0} \sum_{i=0}^{n-1} f(t_i)(g(t_{i+1}) - g(t_i))$$

*converges to $\int_0^\cdot f\,d^y g$ when the mesh $\delta$ of the partition*

$$\pi = \{0 = t_0 < t_1 < \cdots < t_n = t\}$$

*goes to zero.*

Proposition 4.28 permits us to identify the Young integral and the integral of [31]; see Theorem 4.2.1. We thus are allowed to use the following Itô formula established in [31], Theorem 4.3.1, page 351.

PROPOSITION 4.29. *Let $f$ be in $C^\gamma$ and $F$ be in $C^1([0,1] \times \mathbb{R})$ such that $t \mapsto \partial_x F(t, f(t))$ belongs to $C^\alpha$ with $\alpha + \gamma > 1$. Then*

$$F(t, f_t) = F(0, f_0) + \int_0^t \partial_x F(s, f_s)\,d^\circ f_s + \int_0^t \partial_s F(s, f_s)\,ds.$$

We will need the hypothesis

(H$_1'$)     $\sigma$ is in $C^{1,0}$;
        $|\sigma(t,x) - \sigma(t,y)| \le c_n|x - y|^\delta$      $\forall t \in [0,1], |x| + |y| \le n,$

for every $n$ in $\mathbb{N}$, with $c, c_n > 0$, $\delta > \frac{1}{\gamma} - 1$.

We state the proposition, in the Hölder case, which is equivalent to Proposition 4.8, in the finite cubic variation case.

PROPOSITION 4.30. *Let $\sigma$ satisfy (H$_1$), (H$_1'$) and (H$_2$). Suppose that either $P(\eta \in S) = 1$ or $\alpha$ and $\sigma$ are autonomous. Then equation (34) has a unique solution with $\gamma$-Hölder continuous paths, if and only if the following stochastic differential equation has a unique solution:*

(35)          $$Y = \nu^\sigma + \xi + \int_0^\cdot (\widetilde{\partial_s H} + \widetilde\alpha)(s, Y_s)\,ds.$$

We observe that since $\xi$ is $\gamma$-Hölder with $\gamma$ greater than $\frac{1}{2}$, its cubic variation is equal to zero, then equation (35) agrees with equation (26).

REMARK 4.31. Hypothesis (H$_2$) on the the zeros of $\sigma$ is indeed necessary for uniqueness. Suppose $\alpha = 0$, $\sigma$ autonomous and vanishing only at some point $x_0$ with $\frac{1}{\sigma}$ being integrable around $x_0$. Then problem

$$\begin{cases} d^\circ X_t = \sigma(X_t)\,d^\circ \xi_t, \\ X_0 = x_0, \end{cases}$$

has at east two solutions $X_t^1 \equiv x_0$ and $X_t^2 = K(\xi_t)$, where $K = H^{-1}$ and $H(x) = \int_{x_0}^x \frac{1}{\sigma(z)}\,dz$.

COROLLARY 4.32. *Suppose that in addition to the assumptions of Proposition 4.30, $\alpha$ is bounded and locally Lipschitz in $x$ uniformly in $t$, and that $\sigma$ verifies*

$$\sup_{(t,x)\in[0,1]\times S^n} |\partial_t \log(|\sigma(t,x)|)| \leq a_n$$

*for some sequence of positive number $(a_n)_{n\in\mathbb{N}}$. Then equation (34) has a unique solution.*

4.8. *The case of the fractional Brownian motion.* In this section we investigate a significant particular case. We suppose that $\xi = (B_t^H, 0 \leq t \leq 1)$ is a *fractional Brownian motion* on the given filtered probability space (the filtration $\mathbb{F}$ being generated by $B^H$ and the sets of zero probability), with Hurst parameter $H$ strictly larger than $\frac{1}{2}$. Furthermore, we assume that $\eta$ is deterministic, and $\alpha \colon [0,1] \times \mathbb{R} \to \mathbb{R}$ is measurable and locally bounded in $x$, uniformly in $t$. It is well known that $B^H$ has $\lambda$-Hölder continuous paths, for every $\lambda < H$, on $[0,1]$, almost surely. The information about the law $B^H$ allows us to make use of some recent results about uniqueness and existence of a stochastic differential equation driven by a *fractional Brownian motion* with drift equal to 1, which can be found in [22]. More precisely, there the authors establish existence and uniqueness of the integral equation

$$Y_t = y + B_t^H + \int_0^t b(s, Y_s)\, ds, \qquad 0 \leq t \leq 1, y \in \mathbb{R},$$

under this regularity assumption on $b$:

(36)     (H$_4'$)     $|b(t,y) - b(s,x)| \leq C(|x-y|^\alpha + |t-s|^\beta)$,

for some positive constant $C$, with $1 > \alpha > 1 - \frac{1}{2H}$, $\beta > H - \frac{1}{2}$.

Imposing conditions ensuring that the assumption above is satisfied by the coefficients of equation (35), we get the following corollary.

COROLLARY 4.33. *Let $\sigma$ be bounded, satisfying assumptions (H$_1$) and (H$_2$).*

1. *If $P(\eta \in D) = 1$ and both $\alpha$ and $\sigma$ are autonomous, then the integral equation*

(37)     $$X = \eta + \int_0^{\cdot} \sigma(s, X_s)\, d^\circ B_s^H + \int_0^{\cdot} \sigma\alpha(s, X_s)\, ds$$

*has the unique solution $X \equiv \eta$.*

2. *Suppose that $P(\eta \in S^n) = 1$, for some $n$ in $\mathbb{N}$. Let $\alpha$ satisfy hypothesis (H$_4'$), $\sigma$ hypothesis (H$_1'$) and:*

(H$_3'$)
   (i)  $\int_{S^n} |g(t,x) - g(s,x)|\, dx \leq a_n |t-s|^\beta$,
   (ii) $\int_x^y \sup_{t\in[0,1]} |g(t,z)|\, dz \leq a_n |x-y|^\alpha, x, y \in S^n, x \leq y$
   (iii) $\int_{S^n} \sup_{t\in[0,1]} |g(t,z)|\, dz < +\infty$,

*with*

$$g(t,x) = \frac{\partial_t \sigma(t,x)}{(\sigma(t,x))^2}, \qquad (t,x) \in [0,1] \times S,$$

*for some positive constant $a_n$. Then the integral equation (37) has a unique solution.*

Proof. Suppose $\eta \in S^n$. Condition (iii) of (H$_3'$) and the boundness of $\sigma$ imply that $(t,y) \mapsto K^n(t,y)$ is Lipschitz in $x$, uniformly in $t$, and Lipschitz in $t$ uniformly in $x$. Thanks to conditions (i) and (ii), $(t,x) \mapsto \partial_t H^n(t,x)$ fulfills assumption (H$_4'$) for some $C$ positive constant. Then equation (35) has a unique solution by the mentioned result of [22]. Proposition 4.30 permits us to conclude. If $\eta \in D$, uniqueness follows by Proposition 4.30. □

4.9. *Existence in the case of Brownian motion.* If $H = \frac{1}{2}$ and $B^H = B$ is a Brownian motion, supposing $\sigma$ only continuous, it is possible to find a solution to equation

$$(38) \qquad \begin{cases} d^\circ X_t = \sigma(t,X_t)[d^\circ B + \alpha(t,X_t)\,dt], \\ X_0 = \eta. \end{cases}$$

This can be done using the Itô formula permitting to expand $C^1$ functions of *reversible* semimartingales proved in [26]. We recall the result established by [26] (see also [12]) in the case of Brownian motion.

Definition 4.34. A semimartingale $X$ is a *reversible semimartingale* if the process $\hat{X} = (X_{1-t}, 0 \le t \le 1)$ is a semimartingale.

Proposition 4.35. *Let $X = (X^1, \ldots, X^d)$ be a vector of continuous reversible semimartingales and $f$ in $C^1(\mathbb{R}^d)$. Then*

$$f(X_t) = f(X_0) + \sum_{i=1}^d \int_0^t \partial_i f(X_s)\,d^\circ X_s^i.$$

Then we can state the following.

Proposition 4.36. *Let $\sigma$ satisfy (H$_1$), (H$_2$), $\alpha : [0,1] \times \mathbb{R} \to \mathbb{R}$ be measurable and bounded, and $\eta$ deterministic. Suppose that for every $n$ in $\mathbb{N}$, if $\eta$ is in $S^n$,*

$$\sup_{(t,x) \in [0,1] \times S^n} |\partial_t H^n(t,x)| < +\infty.$$

*Then equation (38) has a solution.*

PROOF.    If $\eta \in D$, $X_t \equiv \eta$ is a solution. Suppose $\eta \in S^n$, for some $n$ in $\mathbb{N}$. Equation

$$Y = H^n(0, \eta) + B + \int_0^{\cdot} (\alpha(t, K^n(s, Y_s)) - \partial_s H(s, K^n(s, Y_s)))\, ds$$

admits a solution since the function $(t, y) \mapsto \alpha(t, K^n(t, y)) - \partial_t H^n(t, K^n(t, y))$ is measurable and bounded; see Theorem 35 of [24]. Using the Girsanov theorem, we find that $Y$ is a Brownian motion under a probability measure $P^*$ equivalent to $P$. Therefore, $Y$ is a reversible semimartingale; see example on page 3 of [26]. Then the Itô formula for reversible semimartingales provides a solution to equation (38):

$$X = K^n(\cdot, Y) = \eta + \int_0^{\cdot} \sigma(t, X_t)\, d^\circ B_t + \int_0^{\cdot} \sigma\alpha(t, X_t)\, dt. \qquad \square$$

REMARK 4.37.    We remark that for such a solution $X$, $\int_0^{\cdot} \sigma(s, X_s)\, d^\circ B_s$ is not a proper Stratonovich integral since $\sigma(\cdot, X)$ may not be a semimartingale.

SCUOLA NORMALE SUPERIORE DI PISA
PIAZZA DEI CAVALIERI 7
I-56126 PISA
ITALY
AND
UNIVERSITÉ PARIS 13
INSTITUT GALILÉE
MATHÉMATIQUES
99, AVENUE J.B. CLÉMENT
F-93430 VILLETANEUSE
FRANCE
E-MAIL: coviellor@msn.com

UNIVERSITÉ PARIS 13
INSTITUT GALILÉE
MATHÉMATIQUES
99, AVENUE J.B. CLÉMENT
F-93430 VILLETANEUSE
FRANCE
E-MAIL: russo@math.univ-paris13.fr
URL: http://www.math.univ-paris13.fr/˜russo